\documentclass[12pt]{amsart}
\usepackage{amssymb,amsmath,amsfonts,amsthm,nomencl,mathrsfs} 
\usepackage[arrow, matrix, curve]{xy}
\usepackage[latin1]{inputenc}
\usepackage{a4wide}
\usepackage{xcolor}

\newcommand{\IC}{\mathbb{C}}
\newcommand{\IR}{\mathbb{R}}

\DeclareRobustCommand\abs[1]{\left\vert #1 \right\vert}

\def\bf{\mathbf}
\def\bb{\mathbb}

\newcommand{\question}[1]{\leavevmode{\marginpar{\tiny%
$\hbox to 0mm{\hspace*{-0.5mm}$\leftarrow$\hss}%
\vcenter{\vrule depth 0.1mm height 0.1mm width \the\marginparwidth}%
\hbox to 0mm{\hss$\rightarrow$\hspace*{-0.5mm}}$\\\relax\raggedright #1}}}

\newcommand{\IEE}{\mathscr{E}}

\newcommand{\IHH}{\mathscr{H}}

\newcommand{\IFF}{\mathscr{F}}
\newcommand{\IAA}{\mathscr{A}}

\newcommand{\dom}{\mathrm{Dom}}

\newcommand{\IP}{\mathbb{P}}

\newcommand{\transport}{\slash\slash}

\newcommand{\Id}{\mathrm{d}}

\newcommand{\f}{\frac}

\newcommand{\R}{\mathbb{R}}
\newcommand{\N}{\mathbb{N}}

\newcommand\newdot{{\kern.8pt\cdot\kern.8pt}}
\def\nbull{{\raise1.5pt\hbox{\bf .}}}
\def\bE{\bb E}

\def\transport{/\hspace*{-3pt}/}

\theoremstyle{plain}            
\newtheorem{theorem}{theorem}[section]
\newtheorem{Lemma}[theorem]{Lemma}
\newtheorem{Corollary}[theorem]{Corollary}
\newtheorem{Theorem}[theorem]{Theorem}
\newtheorem{Proposition}[theorem]{Proposition}

\newtheorem{Conjecture}[theorem]{Conjecture}

\theoremstyle{definition}       
\newtheorem{Definition}[theorem]{Definition}

\numberwithin{equation}{section}

\allowdisplaybreaks[1]

\newcommand{\Hmm}[1]{\leavevmode{\marginpar{\tiny%
			$\hbox to 0mm{\hspace*{-0.5mm}$\leftarrow$\hss}%
			\vcenter{\vrule depth 0.1mm height 0.1mm width \the\marginparwidth}%
			\hbox to 0mm{\hss$\rightarrow$\hspace*{-0.5mm}}$\\\relax\raggedright #1}}}

\begin{document}

\begin{titlepage}

\title[Estimates for the covariant derivative of the heat semigroup on differential forms]{Estimates for the covariant derivative of the heat semigroup on differential forms, and covariant Riesz transforms}

\author[R. Baumgarth]{Robert Baumgarth}
\address{Institut für Mathematik, Universität Leipzig, Augustusplatz 10, 04109 Leipzig, Germany}
\email{Robert.Baumgarth@math.uni-leipzig.de}

\author[B. Devyver]{Baptiste Devyver}
\address{Institut Fourier, Universit\'e Grenoble Alpes, 100 rue des maths 38610 Gières, France}
\email{baptiste.devyver@univ-grenoble-alpes.fr}

\author[B. G\"uneysu]{Batu G\"uneysu}
\address{Institut für Mathematik, Universität Potsdam, Campus Golm, Haus 9, Karl-Liebknecht-Straße 24-25, 14476 Potsdam, Germany}
\email{gueneysu@uni-potsdam.de}

\end{titlepage}

\maketitle 

\begin{abstract}
With $\vec{\Delta}_j\geq 0$ is the uniquely determined self-adjoint realization of the Laplace operator acting on $j$-forms on a geodesically complete Riemannian manifold $M$ and $\nabla$ the Levi-Civita covariant derivative, we prove amongst other things
\begin{itemize}
\item a Li-Yau type heat kernel bound for $\nabla  \mathrm{e}^{ -t\vec{\Delta}_j  }$, if the curvature tensor of $M$ and its covariant derivative are bounded,

\item an exponentially weighted $L^p$ bound for the heat kernel of $\nabla  \mathrm{e}^{ -t\vec{\Delta}_j  }$, if the curvature tensor of $M$ and its covariant derivative are bounded,

\item that $\nabla  \mathrm{e}^{ -t\vec{\Delta}_j  }$ is bounded in $L^p$ for all $1\leq p<\infty$, if the curvature tensor of $M$ and its covariant derivative are bounded,

\item a second order Davies-Gaffney estimate (in terms of $\nabla$ and $\vec{\Delta}_j$) for $\mathrm{e}^{ -t\vec{\Delta}_j  }$ for small times, if the $j$-th degree Bochner-Lichnerowicz potential $V_j=\vec{\Delta}_j-\nabla^{\dagger}\nabla$ of $M$ is bounded from below (where $V_1=\mathrm{Ric}$), which is shown to fail for large times if $V_j$ is bounded.
\end{itemize} 
Based on these results, we formulate a conjecture on the boundedness of the covariant local Riesz-transform $\nabla (\vec{\Delta}_j+\kappa)^{-1/2}$ in $L^p$ for all $1\leq  p<\infty$ (which we prove for $1\leq p\leq 2$), and explain its implications to geometric analysis, such as the $L^p$-Calder\'on-Zygmund inequality. Our main technical tool is a Bismut derivative formula for $\nabla  \mathrm{e}^{ -t\vec{\Delta}_j  }$. 
\end{abstract}

\setcounter{page}{1}

\section{Introduction}

Let $M$ be a smooth connected geodesically complete Riemannian $m$-manifold\footnote{All manifolds are understood to be without boundary, unless otherwise stated.}. The geodesic distance will be denoted by $\varrho(x,y)$ and the induced open balls with $B(x,r)$. Given a smooth vector bundle $\IEE\to M$ carrying a canonically given metric and a canonically given covariant derivative, we denote its fiberwise metric simply by $(\cdot,\cdot)$, with $|\cdot|=\sqrt{(\cdot,\cdot)}$ and its covariant derivative simply with
$$
\nabla: \Gamma_{C^{\infty}}(M,\IEE)\longrightarrow \Gamma_{C^{\infty}}(M,\IEE\otimes T^*M).
$$
These remarks apply in particular to $TM\to M$, $T^*M\to M$ or tensor products thereof. We equip $M$ with the Riemannian volume measure $\Id\mu$; sometimes we will use the following {\em local volume doubling} property for the measure $\Id\mu$, which writes: there exists $C>0$ such that for all $0<r\leq R<+\infty$ and all $z\in M$,

\begin{equation}\label{eq:LVD}\tag{LVD}
\frac{\mu(B(z,R))}{\mu(B(z,r))}\leq C \mathrm{e}^{C R} (R/r)^m.
\end{equation}
By the Bishop-Gromov comparison theorem and the well-known formula for the volume of balls in the hyperbolic space, \eqref{eq:LVD} holds if $\mathrm{Ric}\geq -A^2$ for some $A\geq 0$ (and then the constant $C$ in \eqref{eq:LVD} only depends on $m$ and $A$). A well-known consequence of \eqref{eq:LVD} is the following volume comparison inequality: there is a constant $C>0$ such that for all $t>0$, $x_1,x_2\in M$ and $\varepsilon>0$,

\begin{equation}\label{eq:VC}\tag{$\mathrm{VC}_\epsilon$}
\frac{\mu(B(x_2,\sqrt{t})}{\mu(B(x_1,\sqrt{t}))}\leq C \mathrm{e}^{\frac{Ct}{{\epsilon}}} \mathrm{e}^{\epsilon\frac{\varrho(x_1,x_2)^2}{t}}.
\end{equation}
Indeed, letting $r=\varrho(x_1,x_2)$,

\begin{eqnarray*}\label{eq:VC2}
\frac{\mu(B(x_2,\sqrt{t})}{\mu(B(x_1,\sqrt{t}))} &\leq & \frac{\mu(B(x_1,r+\sqrt{t}))}{\mu(B(x_1,\sqrt{t}))}\\
&\leq & C\left(\frac{r}{\sqrt{t}}+1\right)^m \mathrm{e}^{C(r+\sqrt{t})}.
\end{eqnarray*}
Upon using the elementary inequalities 
$$
\mathrm{e}^Cr\leq \mathrm{e}^{\frac{Ct}{8\epsilon}}\mathrm{e}^{\frac{2\epsilon r^2}{t}},\quad \mathrm{e}^{C\sqrt{t}}\leq C'\mathrm{e}^{Ct}, 
$$
one easily gets \eqref{eq:VC} (with a different value of the constant $C$).

Given a smooth metric vector bundle $\IEE\to M$ we define the Banach spaces $\Gamma_{L^p}(M,\IEE)$ given by equivalence classes of Borel sections $\psi$ in $\IEE\to M$ such that 
$$
\left\|\psi\right\|_{p}:=\left\|\>|\psi|\>\right\|_{p}<\infty,
$$
where $\left\||\psi|\right\|_{p}$ denotes the norm of the function $|\psi|$ with respect to $L^p(M)$. Then $\Gamma_{L^2}(M,\IEE)$ canonically becomes a Hilbert space with scalar product
\begin{align}\label{aspa}
\langle \psi_1,\psi_2\rangle= \int(\psi_1,\psi_2) \Id\mu.
\end{align}
In particular, if $\mathrm{Riem}$ denotes the Riemann curvature tensor, seen for instance as a $(0,4)$-tensor, then we can consider $||\mathrm{Riem}||_\infty$ to be the $||\cdot||_\infty$ norm of 

$$\mathrm{Riem}\in \Gamma_{C^\infty}(M,T^0_4M),$$
where $T^{p}_qM\rightarrow M$ is by definition the vector bundle of tensors of type $(p,q)$. Likewise, seeing $\nabla \mathrm{Riem}$ as a $(0,1+4)$-tensor, we can consider $||\nabla \mathrm{Riem}||_\infty$.
 
Given another smooth metric bundle $\IFF\to M$, the operator norm of a linear map
$$
A: \Gamma_{L^p}(M,\IEE)\longrightarrow \Gamma_{L^q}(M,\IFF)
$$
will be denoted by 
$$
\left\|A\right\|_{p,q}= \sup\left\{\left\|Af\right\|_q:\left\|f\right\|_q\leq 1\right\}\in [0,\infty].
$$
Given a smooth linear partial differential operator
$$
P: \Gamma_{C^{\infty}}(M,\IEE)\longrightarrow \Gamma_{C^{\infty}}(M,\IFF),
$$
its formal adjoint with respect to the scalar products $\langle \cdot,\cdot\rangle$ is denoted by
$$
P^{\dagger}: \Gamma_{C^{\infty}}(M,\IFF)\longrightarrow \Gamma_{C^{\infty}}(M,\IEE).
$$
Let
$$
C^{\infty}(M)\stackrel{\Id_0}{\longrightarrow} \Gamma_{C^{\infty}}(M,T^*M) \cdots\stackrel{\Id_{j-1}}{\longrightarrow}\Gamma_{C^{\infty}}(M,\Lambda^{j}T^*M) \stackrel{\Id_j}{\longrightarrow}\Gamma_{C^{\infty}}(M,\Lambda^{j+1}T^*M)\cdots
$$
be the exterior differential. Then one defines the Laplace-Beltrami operator acting on $0$-forms, and the Hodge Laplacian acting on $j$-forms, respectively, by
\begin{align*}
&\Delta:=\Delta_0 := \Id^{\dagger}_0\Id_0:C^{\infty}(M)\longrightarrow C^{\infty}(M),\\
&\vec{\Delta}_j:= \Id^{\dagger}_j\Id_j+\Id_{j-1}\Id^{\dagger}_{j-1}:\Gamma_{C^{\infty}}(M,\Lambda^jT^*M)\longrightarrow \Gamma_{C^{\infty}}(M,\Lambda^jT^*M).
\end{align*}
The induced direct-sum data will be denoted by
$$
\vec{\Delta}, \Id^{\dagger},\Id_j:\Gamma_{C^{\infty}}(M,\Lambda^jT^*M)\longrightarrow \Gamma_{C^{\infty}}(M,\Lambda^jT^*M),
$$
so that $\vec{\Delta}=(\Id+\Id^{\dagger})^2$.

Note the commutation rules $\Id_j \vec{\Delta}_j =\vec{\Delta}_j \Id_j$ and $\Id^{\dagger}_{j-1} \vec{\Delta}_j =\vec{\Delta}_{j-1} \Id^{\dagger}_{j-1}$. In the case $j=m$ and if $M$ is oriented, the Hodge Laplacian $\vec{\Delta}_m$ is just the conjugate of the scalar Laplacian $\Delta=\Delta_0$ by the Hodge star operator 
$$
\star : \Lambda^0T^*M\longrightarrow \Lambda^m T^*M.
$$
The Bochner-Lichnerowicz formula for the Hodge Laplacian writes
$$
\vec{\Delta}_j = \nabla^{\dagger} \nabla+V_j,
$$
where 
$$
V_j\in \Gamma_{C^{\infty}}(M,\mathrm{End}(\Lambda^jT^*M))
$$
is a fiberwise self-adjoint $0$-th order operator, which satisfies 
\begin{align}\label{potential}
|V_j|\leq C|\mathrm{Riem}|,\quad\text{where $C=C(m)>0$ is a constant that only depends on $m$. }
\end{align}
It is worth noting that for $j=1$ one has $V_1=\mathrm{Ric}^{\mathrm{tr}}$, where the Ricci curvature is read as a section 
$$
\mathrm{Ric}\in\Gamma_{C^{\infty}}(M,\mathrm{End}(TM)),
$$
and its transpose
$$
\mathrm{Ric}^{\mathrm{tr}}\in\Gamma_{C^{\infty}}(M,\mathrm{End}(T^*M))
$$
is defined by duality using the Riemannian metric. \vspace{1mm}

Given a Borel function $\Psi:\IR\to \IR$ and a self-adjoint operator $A$, then $\Psi(A)$ is the self-adjoint operator which is induced by the spectral calculus. Its domain of definition is then denoted by $\dom(\Psi(A))$. As $M$ is geodesically complete, for any $j\in\{0,\cdots,n\}$, $\vec{\Delta}_j$ is essentially self-adjoint \cite{chernoff} in $\Gamma_{L^2}(M,\Lambda^jT^*M)$ when initially defined on $\Gamma_{C^{\infty}_c}(M,\Lambda^jT^*M)$. By a usual abuse of notation, the corresponding self-adjoint realizations will be denoted by $\Delta\geq 0$, resp., $\vec{\Delta}_j\geq 0$ again. By local parabolic regularity, for all square-integrable $j$-forms $\alpha\in \Gamma_{L^2}(M,\Lambda^jT^*M)$, the time dependent $1$-form 
$$
(0,\infty)\times M\ni (t,x)\longmapsto \mathrm{e}^{-t\vec{\Delta}_j  } \alpha(x) \in \Lambda^jT^*_xM
$$
has a smooth representative, which extends smoothly to $[0,\infty)\times M$, if $\alpha$ is smooth. Moreover, there exists a uniquely determined smooth map
$$
(0,\infty)\times M\times M\ni (t,x,y)\longmapsto \mathrm{e}^{-t\vec{\Delta}_j  }(x,y)\in \mathrm{Hom}(\Lambda^jT^*_yM,\Lambda^jT^*_xM)\subset (\Lambda^jT^*M)^*\boxtimes \Lambda^jT^*M, 
$$
\emph{the heat kernel} of $\vec{\Delta}_j$, such that for all $\alpha$ as above, $t>0$, $x\in M$ one has
$$
\mathrm{e}^{-t\vec{\Delta}_j  } \alpha(x)=\int_M \mathrm{e}^{-t\vec{\Delta}_j  }(x,y) \alpha(y) \Id\mu(y). 
$$
Note in particular that
$$
\int_M |\mathrm{e}^{-t\vec{\Delta}_j  }(x,y)|^2 \Id\mu(y)<\infty,\quad \mathrm{e}^{-t\vec{\Delta}_j  }(y,x) = \mathrm{e}^{-t\vec{\Delta}_j  }(x,y)^{\dagger},
$$
and that the heat kernel satisfies the usual semigroup identity (cf. Theorem II.1 in \cite{batu}).

Estimates for the heat kernel of the Hodge Laplacian and its derivatives under curvature assumptions, or more generally estimates of the heat kernel of covariant Schrödinger operators of the form $H_V:=\nabla^\dagger\nabla+V$ in $\Gamma_{L^2}(M,\IEE)$, with a smooth metric vector bundle $\IEE\to M$, have already been studied for a long time. When $H_V$ is bounded from below in the sense of quadratic forms, then $H_V$ is essentially self-adjoint \cite{milatovic}; moreover, if the potential $V$, a pointwise self-adjoint smooth section of $\mathrm{End}(\IEE)\to M$, satisfies $V\geq -a^2$ for some constant $a\in\IR$ (meaning that all $x\in M$ all eigenvalues of $V(x):\IEE_x\to\IEE_x$ are bounded from below by $-a^2$), then semigroup domination \cite{Ber} states that for every $\alpha\in \Gamma_{L^2}(M,\IEE)$, $t>0$, one has
$$
|\mathrm{e}^{-tH_V}\alpha|\leq \mathrm{e}^{a^2t} \mathrm{e}^{-t\Delta}|\alpha|.
$$

As, by the Li-Yau heat kernel estimate \cite{LY}, the assumption $\mathrm{Ric}\geq -A^2$ for some constant $A\in \R$ implies the existence of constants $C_j=C_j(A,m)>0$, $D=D(A,m)>0$ (where $C_2=0$, if $A=0$), such that for all $x,y\in M$, $t>0$, one has
$$
\mathrm{e}^{-t\Delta}(x,y)\leq C_1\mu(B(x,\sqrt{t}))^{-1}\mathrm{e}^{C_2t}\mathrm{e}^{-D\frac{\varrho(x,y)^2}{t}},
$$
in this case semigroup domination implies 
$$
|\mathrm{e}^{-tH_V}(x,y)|\leq C_1\mu(B(x,\sqrt{t}))^{-1}\mathrm{e}^{(a^2+C_2)t}\mathrm{e}^{-D\frac{\varrho(x,y)^2}{t}}.
$$
In particular, assuming $||\mathrm{Riem}||_\infty\leq A$ for some $A>0$ (cf. (\ref{potential})) and using the above result for 
$$
H_V=\vec{\Delta}_j=\nabla^{\dagger}\nabla+V_j,
$$
one gets for every $j=0,\cdots,m$,
\begin{equation}\label{eq:vUE}\tag{\textbf{UE}}
|\mathrm{e}^{-t\vec{\Delta}_j}(x,y)|\leq C\mu(B(x,\sqrt{t}))^{-1}\mathrm{e}^{Ct}\mathrm{e}^{-D\frac{\varrho(x,y)^2}{t}},
\end{equation}
where $C,D>0$ depend only on $A$ and $m$. The commutation rules $\Id_j\vec{\Delta}_j=\vec{\Delta}_{j+1}\Id_j$ and $\Id_j^\dagger\vec{\Delta}_j=\vec{\Delta}_{j+1}\Id_j^\dagger$ can then be used in order to prove that similar pointwise estimates hold for the kernels of $\Id_j\mathrm{e}^{-t\vec{\Delta}_j}$ and $\Id_j^\dagger \mathrm{e}^{-t\vec{\Delta}_j}$:

\begin{Proposition}\label{pro:gradUE}

Assume that there is a constant $A>0$ such that $||\mathrm{Riem}||_\infty\leq A$. Then there exist constants $C=C(A,m)>0$, $D=D(A,m)>0$, such that for all $j\in\{1,\dots,m\}$, $x,y\in M$, $t>0$ one has
\begin{equation}\label{eq:dj}\tag{$\Id$\textbf{UE}}
|\Id_j \mathrm{e}^{-t\vec{\Delta}_j}(x,y)|\leq C\mu(B(x,\sqrt{t}))^{-1}t^{-1/2}\mathrm{e}^{Ct}\mathrm{e}^{-D\frac{\varrho(x,y)^2}{t}},
\end{equation}

\begin{equation}\label{eq:d*j}\tag{$\Id^\dagger$\textbf{UE}}
|\Id_{j-1}^\dagger \mathrm{e}^{-t\vec{\Delta}_j}(x,y)|\leq C\mu(B(x,\sqrt{t}))^{-1}t^{-1/2}\mathrm{e}^{Ct}\mathrm{e}^{-D\frac{\varrho(x,y)^2}{t}}.
\end{equation}
\end{Proposition}

Above and in the sequel, we understand $\Id_j$ to act on the first variable of the heat kernel, so
$$
\Id_j \mathrm{e}^{-t\vec{\Delta}_j}(x,y):= \Id_j \mathrm{e}^{-t\vec{\Delta}_j}(\bullet,y)(x),
$$
and likewise for $\Id^{\dagger}_{j-1}$, and in similar situations such as $\nabla \mathrm{e}^{-t\vec{\Delta}_j}(x,y)$. Note here that for fixed $y$, the map $x\mapsto \mathrm{e}^{-t\vec{\Delta}_j}(x,y)$ becomes a section of a bundle of the form 
$$
(\Lambda^jT^*M)^*\otimes W\cong \Lambda^jT^*M\otimes W\longrightarrow M,
$$
with $W$ a fixed finite dimensional linear space, which explains the action of these differential operators on the heat kernel. 

Although we expect this result to be well-known to the experts, for the sake of completeness, we will provide a proof of Proposition \ref{pro:gradUE} in the appendix.

In this article, our main goal is to prove the analogous estimates for the covariant derivative of the heat kernel of the Hodge Laplacian: more precisely, we wish to obtain pointwise estimates of the form:

\begin{equation}\label{eq:CUE}\tag{$\nabla$\textbf{UE}}
\left|\nabla \mathrm{e}^{-t\vec{\Delta}_j}(x,y)\right|\leq C\mu(B(x,\sqrt{t}))^{-1}t^{-1/2}\mathrm{e}^{Ct}\mathrm{e}^{-D\frac{\varrho(x,y)^2}{t}}.
\end{equation}

For $j=0$ we have $\nabla=d_0$, so \eqref{eq:dj} and \eqref{eq:CUE} are equivalent, and can thus be obtained with assuming merely that $||\mathrm{Riem}||_\infty<\infty$. The same is true for $j=m$ by Hodge duality (if $M$ is oriented). However, for $j\in \{1,\cdots,m-1\}$, the corresponding covariant derivative estimates are significantly stronger than \eqref{eq:dj} and \eqref{eq:d*j}, and are harder to prove as well, as we shall see. In fact, in order to prove these, we will not only need a uniform bound on the Riemannian curvature tensor, but also on its covariant derivative. We can now state our main result:
\begin{Theorem}\label{main} Assume
\begin{align}\label{wlp}
\max\big( \left\|\mathrm{Riem}\right\|_{\infty},\left\|\nabla \mathrm{Riem}\right\|_{\infty} \big)\leq A\quad\text{for some constant $A>0$.}
\end{align} 
Then \eqref{eq:CUE} holds; more precisely, there exist constants $C=C(A,m), D=D(A,m)>0$, such that for all $j\in\{1,\dots,m\}$, $t>0$, $x,y \in M$ one has
\begin{equation}\label{eq:main}
\left|\nabla \mathrm{e}^{-t\vec{\Delta}_j}(x,y)\right|\leq C\mu(B(x,\sqrt{t}))^{-1}t^{-1/2}\mathrm{e}^{Ct}\mathrm{e}^{-D\frac{\varrho(x,y)^2}{t}}.
\end{equation}
\end{Theorem}

The proof of Theorem \ref{main} is given in Section \ref{beweis} and is based on a probabilistic representation of $\nabla \mathrm{e}^{-\frac{t}{2} \vec{\Delta}_j}(x,y)$ in terms of the Brownian bridge, namely a so called Bismut derivative formula, which should be of independent interest and which is proved based on the methods from \cite{driver} (see also \cite{wang}) in Section \ref{bdf}. In fact, we first prove a local Bismut derivative formula for $\nabla \mathrm{e}^{-\frac{t}{2} \vec{\Delta}_j}\alpha(x)$ in terms of Brownian motion for $\alpha\in \Gamma_{C^{\infty}\cap L^2}(M,\Lambda^j T^*M)$, which does not require any assumptions on the geometry. Then we use this formula to obtain global $L^{\infty}$ estimates under (\ref{wlp}) for $\nabla \mathrm{e}^{-\frac{t}{2} \vec{\Delta}_j}$, which are then used to prove a global Bismut derivative formula for $\nabla \mathrm{e}^{-\frac{t}{2} \vec{\Delta}_j}\alpha(x)$ in terms of Brownian motion. The reason for this rather technical procedure is that, unlike its global counterpart, the local Bismut derivative formula contains a first exit time of Brownian motion from a ball $B$ around its starting point $x$, which is why this formula cannot be controlled well in terms of Brownian bridge (which is conditioned to be in a fixed point $y$ at its terminal time, which need not be in in $B$). 

At this point, let us mention the recent paper \cite{KP}, where, with completely different methods, $L^p\to L^q$ estimates for the covariant derivates of heat kernels of covariant Schrödinger operators have been considered for so called \emph{asymptotically locally Euclidean Riemannian manifolds}. In addition to the fact a very special form of the geometry is required, the estimates from \cite{KP} are also different in their nature than ours: in \cite{KP} one needs an additional decay at $\infty$ of the potentials in order to obtain a damping effect in the constants also for for large $t>0$, while our estimates do not require any decay at $\infty$ of the potentials, while they only damp for small $t>0$.
\vspace{1mm}

As a consequence of the pointwise estimates from Theorem \ref{main} and local volume doubling, one obtains $L^p\to L^p$ bounds for $\nabla \mathrm{e}^{-t\vec{\Delta}_j}$, as well as weighted $L^p$-estimates for the kernel of $\nabla \mathrm{e}^{-t\vec{\Delta}_j}$:

\begin{Corollary}\label{cor:gradient}
Assume 
\begin{align*}
\max\big( \left\|\mathrm{Riem}\right\|_{\infty},\left\|\nabla \mathrm{Riem}\right\|_{\infty} \big)\leq A\quad\text{for some constant $A>0$.}
\end{align*} 
Then:\\
\emph{I)} For all $1\leq p<\infty$ there exists a constant $C=C(A,m,p)>0$, such that for all $j\in\{1,\dots,m\}$, $t>0$ one has
$$
\left\| \nabla  \mathrm{e}^{- t\vec{\Delta}_j }\right\|_{p,p}\leq C\mathrm{e}^{tC} t^{-1/2}.
$$
\emph{II)} There exists a constant $\gamma=\gamma(A,m)>0$, and for all $1\leq p<\infty$ a constant $C=C(A,m,p)>0$, such that for all $j\in\{1,\dots,m\}$ and $t>0$ one has 
\begin{align}\label{suppo}
\int |\nabla \mathrm{e}^{-t\vec{\Delta}_j}(x,y)|^p \mathrm{e}^{\frac{\gamma\varrho(x,y)^2}{t}}\Id\mu(x)\leq \frac{C\mathrm{e}^{Ct}}{t^{p/2}\mu( B(y,\sqrt{t}) )^{p-1}}.
\end{align}

\end{Corollary}

Corollary \ref{cor:gradient} is proved in Section \ref{beweis2}.\vspace{1mm}

It is a well-known principle that stems from the work of Coulhon and Duong \cite{cd}, as well as later works by Auscher, Coulhon, Duong and Hofmann \cite{ACDX}, that (at least for the scalar Laplacian) estimates for the spatial derivative of the heat kernel should have consequences for the corresponding Riesz transform. In this respect, applying (\ref{suppo}) with $p=2$, we are going to establish the following result concerning the {\em covariant} Riesz transform:

\begin{Corollary}\label{cor:Riesz} Assume
\begin{align*}
\max\big( \left\|\mathrm{Riem}\right\|_{\infty},\left\|\nabla \mathrm{Riem}\right\|_{\infty} \big)\leq A\quad\text{for some constant $A>0$.}
\end{align*} 
Then there exists a $\kappa_0=\kappa_0(A,m)>0$, which only depends on $A$, $m$, such that for all $\kappa\geq \kappa_0$, and all $j\in\{1,\dots,m\}$, the operator $\nabla (\vec{\Delta}_j+\kappa)^{-1/2}$ is of weak $(1,1)$ type with a bound only depending on $A,m,\kappa$; in other words, there exist a constant $D=D(A,m,\kappa)>0$, which only depends on $A$, $m$ and $\kappa$, such that for all $j\in\{1,\dots,m\}$, $\lambda>0$, $f\in\Gamma_{L^1}(M,\Lambda^jT^*M)$ one has 
$$
\mu\{|\nabla (\vec{\Delta}_j+\kappa)^{-1/2}f|>\lambda\}\leq\frac{D}{\lambda} \left\|f\right\|_1.
$$
In particular, for all $1<  p\leq  2$ there exists a constant $C=C(A,m,p,\kappa)>0$, which only depends on $A$, $m$, $p$, $\kappa$, such that for all $j\in\{1,\dots,m\}$ one has
\begin{align}\label{apodddy}\tag{$\nabla {\bf R}_p$}
\left\| \nabla (\vec{\Delta}_j+\kappa)^{-1/2} \right\|_{p,p}\leq C.
\end{align}
\end{Corollary}

Corollary \ref{cor:Riesz} is proved in Section \ref{RE} (where we show that the (1,1) property indeed implies the $L^p$-boundedness through). This results improves a result by Thalmaier and Wang \cite[Theorem D]{wang}: more precisely, in \cite[Theorem D]{wang} the same conclusion for the covariant Riesz transform is obtained, however an additional assumption on the volume growth of $M$ is made. This volume assumption excludes in particular hyperbolic geometries (see \cite{Pig}), while such geometries are covered by our Corollary \ref{cor:Riesz}. In light of the our main result, Theorem \ref{main}, and the results in \cite{ACDX} for the {\em scalar} Riesz transform, it is natural to expect that a uniform bound on $\mathrm{Riem}$ and $\nabla\mathrm{Riem}$ implies that the covariant Riesz transform is bounded on $L^p$ for all $1<p<\infty$; specifically, we make the following conjecture:

\begin{Conjecture}\label{conj:riesz}Assume
\begin{align*}
\max\big( \left\|\mathrm{Riem}\right\|_{\infty},\left\|\nabla \mathrm{Riem}\right\|_{\infty} \big)\leq A<\infty.
\end{align*} 
There exists a $\kappa_0=\kappa_0(A,m)>0$, which only depends on $A$, $m$, such that for all $\kappa\geq \kappa_0$, $j\in \{1,\cdots,m\}$, $p\in (1,\infty)$, $0<\kappa\leq \kappa_0$ one has $\left\|\nabla (\vec{\Delta}_j+\kappa)^{-1/2}\right\|_{p,p}<\infty$, with bound only depending on $A,m,\kappa,p$.
\end{Conjecture}

It should be noted that in \cite{Bakry}, one can find the following result concerning the covariant Riesz transform for {\em Einstein} manifolds:

\begin{Theorem}{Bakry, \cite[Theorem 6.1]{Bakry}} Assume $M$ is Einstein and $||\mathrm{Riem}||_\infty\leq A$ for some constant $A>0$. Then for every $p\in (1,\infty)$, there exists a $\kappa_0=\kappa_0(p,A,m)>0$, which only depends on $p$, $A$, $m$, such that for all $\kappa\geq \kappa_0$, $j\in\{1,\dots, m\}$ one has $\left\|\nabla (\vec{\Delta}_j+\kappa)^{-1/2}\right\|_{p,p}<\infty$, with a norm bound that only depending on $A,p,\kappa,m$.
\end{Theorem}

However, the case of Einstein manifold is very special, because for Einstein manifolds there is a nice commutation formula beween $\nabla$ and $\vec{\Delta}_j$ (cf. \cite[formula (6.1)]{Bakry}). We currently do not know whether the assumption on $\nabla \mathrm{Riem}$ is really necessary in Conjecture \ref{conj:riesz}; however, it is known that the curvature hypotheses cannot be weakened to merely boundedness from below of the sectional curvature: in fact a recent result of Marini and Veronelli \cite{MV} shows that there exist manifolds with positive sectional curvature, for which the covariant Riesz transform is not bounded on $L^p$ for all $p\in (1,\infty)$. 

\medskip

It should also be noted that boundedness in $L^p$ of the Riesz transforms $ \Id_j (\vec{\Delta}_j+\kappa)^{-1/2}$ and $ \Id_{j-1}^\dagger(\vec{\Delta}_j+\kappa)^{-1/2}$ instead of $\nabla (\vec{\Delta}_j+\kappa)^{-1/2}$ is again a considerably easier business. In fact, a classical result by Bakry \cite[Theorem 5.1]{Bakry} states:

\begin{Theorem}\label{thm:riesz-bakry} Assume $||\mathrm{Riem}||_\infty\leq A$ for some constant $A>0$; then there exists a $\kappa_0=\kappa_0(A,m)>0$, which only depends on $A$, $m$, such that for all $\kappa\geq \kappa_0$, and all $j\in \{0,\dots,m\}$ the operators $\Id_j (\vec{\Delta}_j+\kappa)^{-1/2}$ and $\Id_{j-1}^\dagger(\vec{\Delta}_j+\kappa)^{-1/2}$ are weakly $(1,1)$ with an $(1,1)$-norm bound that only depends on $A$, $m$, $\kappa$; in particular, for every $p\in (1,\infty)$, one has
$$
\left\|\Id_j(\vec{\Delta}_j+\kappa)^{-1/2}\right\|_{p,p}<\infty, \qquad \left\|\Id_{j-1}^\dagger(\vec{\Delta}_j+\kappa)^{-1/2}\right\|_{p,p}<\infty,
$$
with norm bounds only depending on $A,p,\kappa,m$.
\end{Theorem}
As we said, the $L^p$-boundedness part of this result follows from \cite[Theorem 5.1]{Bakry}; the weak $(1,1)$ part appears to be new in this generality. The latter is established using the estimates \eqref{eq:dj} and \eqref{eq:d*j}, and Coulhon-Duong theory as in the proof of Corollary \ref{cor:Riesz}, yielding an alternative proof of the $L^p$ boundedness part of Theorem \ref{thm:riesz-bakry}. This will be done in Section \ref{RE}.

\medskip

Let us stress that for applications in geometric analysis, the $L^p$-boundedness of $\nabla (\vec{\Delta}_j+\kappa)^{-1/2}$ is more important than that of $\Id_j (\vec{\Delta}_j+\kappa)^{-1/2}$ or $\Id_j^\dagger(\vec{\Delta}_j+\kappa)^{-1/2}$. For example, as shown in \cite[Proof of Theorem 4.13]{gp1}, the former boundedness for $j=1$ implies the $L^p$-Calder\'on-Zygmund inequality
\begin{equation}\label{eq:CZ}\tag{$\mathrm{CZ}_p$}
\left\| \mathrm{Hess} ( u)\right\|_p\leq D_{\mathrm{CZ}} (\left\| \Delta  u\right\|_p+ \left\|   u\right\|_p)\quad\text{ for all $u\in C^{\infty}_c(M)$},
\end{equation}
where $D_{\mathrm{CZ}}$ only depends on $\left\|\nabla ( \vec{\Delta}_j+\kappa)^{-1/2}\right\|_{p,p}$.\vspace{1mm}

The $L^p$-Calder\'on-Zygmund inequality together with $\left\|\mathrm{Riem}\right\|_{\infty}<\infty$, in turn, implies global a-priori $L^p$-estimates (cf. Theorem 4 b) in \cite{gp2}) for distributional solutions $\Psi\in L^p(M)$ of the Poisson equation $\Delta \Psi=f\in L^p(M)$ which is of the form
$$
\left\| \mathrm{Hess} ( \Psi) \right\|_p+\left\|  \nabla  \Psi\right\|_p\leq   C(\left\|f\right\|_p+\left\|\Psi\right\|_p),
$$
where $C$ only depends on $D_{\mathrm{CZ}}$ and any upper bound for $\left\|\mathrm{Riem}\right\|_{\infty}$. Hence, Corollary \ref{cor:Riesz} readily implies:

\begin{Corollary}\label{cz}
Assume
\begin{align*}
\max\big( \left\|\mathrm{Riem}\right\|_{\infty},\left\|\nabla \mathrm{Riem}\right\|_{\infty} \big)\leq A\quad\text{for some constant $A>0$.}
\end{align*} 
Then for all $1< p\leq  2$, there exists a constant $D'=D'(A,m,p)>0$ such that
\begin{align}\label{aappqq}
\left\| \mathrm{Hess} ( u)\right\|_p\leq D' (\left\| \Delta  u\right\|_p+ \left\|   u\right\|_p)\quad\text{for all $u\in C^{\infty}_c(M)$,}
\end{align}
and such that for every distributional solution $\Psi\in L^p(M)$ of $\Delta \Psi=f\in L^p(M)$ one has
$$
\left\| \mathrm{Hess} (\Psi) \right\|_p+\left\| \nabla  \Psi\right\|_p\leq   D'(\left\|f\right\|_p+\left\|\Psi\right\|_p).
$$
\end{Corollary}

Note that the CZ-inequality (\ref{aappqq}) improves Theorem D in \cite{gp1} by getting rid of the volume assumption made there.\vspace{2mm}

In addition to a priori estimates for the Poisson equation, the $L^p$-Calder\'on-Zygmund inequality implies precompactness results for isometric immersions (cf. Theorem 1.1 in \cite{breuning} and \cite{gp2}). Moreover, the recent survey article \cite{Pig} contains the state-of-the art for the $L^p$-Calder\'on-Zygmund inequality for large $p$: it is explained therein that $\left\|\mathrm{Riem}\right\|_{\infty}<\infty$ is enough for the $L^p$-Calder\'on-Zygmund inequality to hold for all $p>\max(2,m/2)$. In this sense, Corollary \ref{cz} can be considered a complementary result for small $p$.\vspace{2mm}

A fundamental tool in \cite{ACDX} for obtaining the boundedness of the Riesz transform on $L^p$ for $p>2$ are the so-called {\em Davies-Gaffney estimates} for the gradient of the scalar Laplacian, that is to say $L^2$ off-diagonal estimates for $\mathrm{e}^{-t\Delta}$ and $\Id\mathrm{e}^{-t\Delta}$. At zeroth order, these estimates are equivalent to the finite speed of propagation of the associated wave equation, and they hold true for $\mathrm{e}^{-t\vec{\Delta}_j}$ for all $j=0,\cdots,m$ \cite{Sik} (note that Davies-Gaffney estimates for covariant Schrödinger semigroups  of the form $\mathrm{e}^{-tH_V}$ for unbounded $V$'s play a fundamental role in the context of essential self-adjointness of covariant Schrödinger operators \cite{batu2}). One can ask more generally whether Davies-Gaffney hold for the covariant derivative of the heat kernel of the Hodge Laplacian. In this respect, we have the following result, which is proved in Section \ref{beweis3}, and where $\chi_A$ denotes the indicator function of a set $A\subset M$:

\begin{Theorem}\label{main2} There exist universal constants $c_1,c_2>0$ such that for all $j\in\{1,\dots,m\}$ with $V_j\geq - A$ for some constant $A\geq 0$, all $t>0$, all Borel subsets $E,F\subset M$ with compact closure, and all $\alpha\in \Gamma_{L^2}(M,\Lambda^jT^*M)$ with $\mathrm{supp}(\alpha)\subset E$, one has
$$
\left\|\chi_F \mathrm{e}^{-t \vec{\Delta}_j}\alpha\right\|_2+\left\|\chi_F\sqrt{t} \nabla \mathrm{e}^{-t \vec{\Delta}_j}\alpha\right\|_2+\left\|\chi_Ft\vec{\Delta}_j  \mathrm{e}^{-t \vec{\Delta}_j}\alpha\right\|_2\leq c_1(1+\sqrt{t}A)\mathrm{e}^{-\frac{c_2 \varrho(E,F)^2}{t}} \left\|\chi_E\alpha\right\|_{2}.
$$ 
\end{Theorem}
Actually, the above Davies-Gaffney estimate for $\mathrm{e}^{-t \vec{\Delta}_j}$ and $t\vec{\Delta}_j\mathrm{e}^{-t \vec{\Delta}_j}$, even without the extra $\sqrt{t}$ factor on the right-hand side, are already known (cf \cite[Lemma 3.8]{AMR}), but for the sake of completeness we will provide a proof. The novelty is the Davies-Gaffney bound for the gradient term $\sqrt{t} \nabla \mathrm{e}^{-t \vec{\Delta}_j}$.

Note also that the above Davies-Gaffney bounds implies that for all $0<t<1$ one has
$$
\left\|\chi_F \mathrm{e}^{-t \vec{\Delta}_j}\alpha\right\|_2+\left\|\chi_F\sqrt{t} \nabla \mathrm{e}^{-t \vec{\Delta}_j}\alpha\right\|_2+\left\|\chi_Ft\vec{\Delta}_j  \mathrm{e}^{-t \vec{\Delta}_j}\alpha\right\|_2\leq c_{1,A}\mathrm{e}^{-\frac{c_2 \varrho(E,F)^2}{t}} \left\|\chi_E\alpha\right\|_{2},
$$ 
which is ultimately what is needed for the machinery from \cite{ACDX}. Remarkably, for $j>0$ the latter form of the Davies inequality is false for large times, unless one makes additional geometric assumptions on $M$. This means that, contrary to what happens for the scalar Laplacian, even $L^2$ off-diagonal estimates for the covariant derivative of the heat operator of the Hodge Laplacian are non-trivial. This is the content of the following result, which is proved in Section \ref{beweis4}: 

\begin{Theorem}\label{main3} Assume that $M$ is noncompact, that there exists $j\in\{1,\dots, m\}$ with $\left\|V_j\right\|_\infty<\infty$, and that there exist constants $c_1,c_2>0$ such that for all $t>0$, all Borel subsets $E,F\subset M$ with compact closure and all $\alpha\in \Gamma_{L^2}(M,\Lambda^jT^*M)$ with $\mathrm{supp}(\alpha)\subset E$, one has
$$
\left\|\chi_F\sqrt{t} \nabla \mathrm{e}^{-t \vec{\Delta}_j}\alpha\right\|_2 \leq c_1\mathrm{e}^{-\frac{c_2 \varrho(E,F)^2}{t}} \left\|\chi_E\alpha\right\|_{2}.
$$
Then one has $\mathrm{Ker}_{L^2}(\vec{\Delta}_j)=\{0\}$.
\end{Theorem}

\medskip
 
 
\medskip

\textbf{Acknowledgements: } The authors are indepted to Stefano Pigola and Anton Thalmaier for very helpful discussions.\\

B. Devyver was partly supported by the French ANR through the project RAGE ANR-18-CE40-0012, and as well as in the framework of the ``Investissements d'avenir'' program (ANR-15-IDEX-02) and the LabEx PERSYVAL (ANR-11-LABX-0025-01).

\section{Bismut derivative formula}\label{bdf}

Fix $j\in\{1,\dots,m\}$. The following endomorphisms are built from the curvature and its first derivative and will play a crucial role in the probabilistic formula for $\nabla  \mathrm{e}^{-\frac{t}{2}\vec{\Delta}_j }$, the main result of this section. In this section, we read the Riemannian curvature as a section
$$
\mathrm{Riem} \in \Gamma_{C^{\infty}} (M,  T^*M\otimes T^*M  \otimes \mathrm{End}(TM) ).
$$
Then the section
$$
\underline{V}_j   \in \Gamma_{C^{\infty}}(M,\mathrm{End}(T^*M\otimes \Lambda^jT^*M))=\Gamma_{C^{\infty}}\big(M,\mathrm{End}\big(\mathrm{Hom}(TM, \Lambda^jT^*M)\big)\big)
$$
is defined on $x\in M$, $\phi\in T_x^*M\otimes \Lambda^jT_x^*M$, $v\in T_x^*M$, by
$$
\underline{V}_j(\phi)(v)= (\mathrm{Ric}^{\mathrm{tr}}\otimes 1_{\Lambda^jT^*_xM})(\phi)(v)+( 1_{T^*_xM}\otimes V_j)(\phi)(v)-2\sum^m_{i=1}\mathrm{Riem}(v,e_i )\phi(v),
$$
where $e_j$ is any smooth local orthonormal basis for $T_x M$, and the section
$$
\rho_j   \in \Gamma_{C^{\infty}}(M,\mathrm{Hom}(T^*M, T^*M\otimes T^*M))
$$
is defined on $\alpha\in \Lambda^jT^*_xM$, $v\in T_xM$ by
	\begin{align*}
    	\rho_j (\alpha)(v) = (\nabla_v V_j)\alpha  +\sum^m_{i=1}(\nabla_{e_i}\mathrm{Riem}^{\mathrm{tr}})(e_i,v)\alpha.
	\end{align*}

For the formulation of the probabilistic results of this section, we will assume that the reader is familiar with stochastic analysis on manifolds. Classical references in this context are e.g. \cite{ikeda,hsu, elworthy,hackenbroch} (see also \cite{boldt}	for a very brief summary the notions relevant in the sequel).\\
Let $(\Omega, \mathcal{F}, \mathcal{F}_*, \IP)$ be a filtered probability space which satisfies the usual conditions and which for every $x\in M$ carries an adapted Brownian motion
$$
X^x: [0,\zeta^x)\times \Omega \longrightarrow M
$$
starting from $x\in M$, where 
$$
\zeta^x:\Omega \longrightarrow (0,\infty]
$$
denotes the lifetime of $X^x$ (noting that $\zeta^x=\infty$ a.s., if for example $\mathrm{Ric}\geq -a$ for some $a>0$). Given a metric vector bundle $E\to M$ with metric connection, let  
$$
\transport^x: [0,\zeta^x)\times \Omega\longrightarrow   \mathrm{Hom}(E_{x} , (X^x)^*E)
$$
denote the (pathwise orthogonal) parallel transport with respect to $\nabla$ along $X^x$. We define continuous adapted processes with paths having a locally finite variation by
\begin{align*}&Q^x_j:[0,\zeta^x)\times \Omega\longrightarrow\mathrm{End}(\Lambda^jT_{x}^*M),\\
(\Id/ \Id s) Q^x_j(s) &= - \frac{1}{2} Q^x_j(s)  \big(\transport_s^{x,-1} V_j(X_s^x)\transport_s^x \big) ,\quad    Q^x_j(0) = 1_{ \Lambda^jT_{x}^* M},
    \end{align*}
and
\begin{align*}
&\underline{Q}^x_j :[0,\zeta^x )\times \Omega\longrightarrow\mathrm{End}(T_{x}^*M\otimes \Lambda^jT_{x}^*M ),\\
(\Id/\Id s) \underline{Q}^x_j(s) &= - \frac{1}{2} \underline{Q}^x_j(s) \big( \transport_s^{x,-1} \underline{V_j}(X_s^x)\transport^x_s \big),\quad    \underline{Q}^x_j(0)  = 1_{ T^*_{x} M\otimes \Lambda^jT^*_{x} M}.
    \end{align*}
 
In addition, for every $r>0$ let
$$
\tau^x_r :=\inf\{t\in [0,\zeta^x): X^x_{t}  \notin B(x,r)\}: \Omega \longrightarrow [0,\infty]
$$
be the first exit time of $X^x $ from $B(x,r)$. Note that $\zeta^x>\tau^x_r >0$ $\mathbb{P}$-a.s. Let
$$
\underline{X}^x:=  \int_0^{\bullet} \transport_s^{x,-1}  \circ \Id X_s^x:[0,\zeta^x)\times \Omega\longrightarrow T_{x}M
$$
denote the anti-development of $X^x$ (a Brownian motion in the Riemannian manifold $T_xM$). Here, $\circ \Id $ denotes the Stratonovich stochastic differential, where It\^{o} stochastic differentials will be denoted by $\Id$. The following definitions will be very convenient in the sequel:

\begin{Definition} Let $r>0$, $t>0$, $x\in M$, $\xi\in T_{x}^*M\otimes \Lambda^jT_{x}M$. We define a set of processes $\IAA_j(x,r,t,\xi)$ to be given by all bounded adapted process 
$$
\ell: [0,t] \times \Omega\longrightarrow T_{x}^*M\otimes \Lambda^jT_{x}M
$$
with locally absolutely continuous paths such that   
$$
\mathbb{E}\left[\int_0^{t\wedge \tau^x_r }\abs{\dot\ell_s}^2 \Id s\right]<\infty,\quad  \ell_0 = \xi,\quad \ell_s = 0\quad \text{for all $s \geq t \wedge\tau^x_r$}.
$$
For every $\ell\in \IAA_j(x,r,t,\xi)$ we define the continuous semimartingale 
\begin{align*}
   & U^{(\ell)} :[0,t\wedge \tau^x_r]\times \Omega\longrightarrow T_{x}M ,\\
		& U^{(\ell)}:=  \int_0^{\bullet} Q^x_j(s)^{\dagger,-1} \Id\underline{X}^x_s \underline{Q}^x_j(s)^{\dagger}  \dot \ell_s+\frac{1}{2} \int_0^{\bullet} Q^x_j(s)^{\dagger,-1}   \big(\transport_s^{x,-1}  \rho_j(X_s^x) \transport_s^x\big)^{\dagger} \underline{Q}^x_j(s)^{\dagger}  \ell_s \Id s.
			\end{align*}
\end{Definition}

The proof of the following result follows the arguments of Theorem 4.1 from \cite{driver}:

\begin{Theorem}[Local covariant Bismut formula]\label{main0} Let $t > 0$, $r>0$, $x\in M$, $\xi\in T_{x}^*M\otimes \Lambda^jT_{x}M$. Then for every $\ell\in\IAA_j(x,r,t,\xi)$ and every $\alpha\in \Gamma_{L^2\cap C^{\infty}}(M,\Lambda^jT^*M)$ one has
\begin{align}\label{eq::nablaestimg}
 (  \nabla  \mathrm{e}^{-\frac{t}{2}\vec{\Delta}_j } \alpha(x)  , \xi   ) = -\bE  \left[ \big( Q^x_j(t\wedge\tau^x_r)  \transport_{t\wedge \tau^x_r }^{x,-1}  \mathrm{e}^{\frac{-(t-t\wedge \tau^x_r)}{2} \vec{\Delta}_j}\alpha(X^x_{t\wedge \tau^x_r } )  , U^{(\ell)}_{t\wedge \tau^x_r } \big)\right].
\end{align}
\end{Theorem}

The following two well-known facts will be used in the proof of Theorem \ref{main0}:

\begin{Lemma}\label{paaq} Let $\tau$ be a $\mathbb{P}$-a.s. finite stopping time and 
$$
Y:[0,\tau]\times \Omega \longrightarrow \IHH
$$ 
be a continuous local martingale taking values in a finite dimensional Hilbert space $\IHH$. Then $Y$ is a martingale, if
$$
\mathbb{E}\left[\sup_{ t\in [0,\tau]} |Y_t|\right]<\infty.
$$
\end{Lemma}

\begin{Lemma}[Burkholder-Davis-Gundy inequality] For all $0<q<\infty$ there exists a constant $C(q)<\infty$ with the following property: if $\tau$ is a $\mathbb{P}$-a.s. finite stopping time and 
$$
Y:[0,\tau]\times \Omega \longrightarrow \IHH
$$ 
is a continuous local martingale taking values in a finite dimensional Hilbert space $\IHH$ and staring from $0$, then one has
$$
\mathbb{E}\Big[(\sup_{ t\in [0,\tau]}|Y_t|)^q\Big]\leq C(q) \mathbb{E}\big[|[Y]_{\tau}|^{q/2}\big],
$$
where 
$$
[Y]:[0,\tau]\times \Omega \longrightarrow \mathrm{End}(\IHH)
$$
denotes the quadratic variation of $Y$ and $|[Y]_{\tau}|$ the Hilbert-Schmidt norm. 
\end{Lemma}

\begin{proof}[Proof of Theorem \ref{main0}] Pick $A>0$ such that
$$
\max\big( \left|\mathrm{Riem}\right|,\left|\nabla \mathrm{Riem}\right| \big)\leq A\quad\text{in $B(x,r)$.}
$$
We start by noting that for all $s\geq 0$ one has
\begin{align} \label{qbit1}
|Q_j^x(s) |=|Q_j^x(s)^{\dagger}|\leq \mathrm{e}^{C (m, A)s}, \quad |\underline{Q}^x_j(s)|=|\underline{Q}^x_j(s)^{\dagger} |\leq  \mathrm{e}^{C (m,A) s}\quad \text{ $\IP$-a.s. in $\{ s \leq \tau^x_r\}$ },
\end{align}
by Gronwall\rq{}s lemma, and as $Q^x_j $ and $\underline{Q}^x_j$ are invertible with
 \begin{align*}
(\Id/\Id s) Q_j^x(s)^{-1}  =  \frac{1}{2}   \big(\transport_s^{x,-1} V_j(X_s^x)\transport_s^x\big)Q_j^x(s)^{-1} ,\quad    Q_j^x(0)^{-1}  = 1_{ \Lambda^jT^*_{x} M},
    \end{align*}
and
\begin{align*}
(\Id/\Id s) \underline{Q}_j^x(s)^{-1} =  \frac{1}{2}  \Big( \transport_s^{x,-1}\underline{V}_j(X^x_s)\transport_s^x\Big) \underline{Q}_j^x(s)^{-1},\quad     \underline{Q}_j^x(0)^{-1} = 1_{ T^*_{x} M\otimes \Lambda^j T^*_{x} M},
\end{align*}
we also have 
\begin{align}\label{qbit2}
|Q_j^x(s)^{-1} |=|Q_j^x(s)^{-1,\dagger} |\leq \mathrm{e}^{C (m ,A)s}, \quad |\underline{Q}_j^x(s)^{-1}|=|\underline{Q}_j^x(s)^{-1,\dagger}|\leq \mathrm{e}^{C (m,A) s}\quad \text{ $\IP$-a.s. in $\{ s \leq \tau^x_r\}$ }.
\end{align}
Using It\^{o}'s formula one shows that
\begin{align*}
   Y &:= \big( Q^x_j   \transport ^{x,-1}  \nabla \mathrm{e}^{-\frac{t-\bullet}{2}\vec{\Delta}_j}  \alpha(X^x ) , \ell\big) - \big( Q^x_j   \transport^{x,-1} \mathrm{e}^{-\frac{t-\bullet}{2}\vec{\Delta}_j} \alpha(X^x ) , U^{(\ell)} \big) \\
   &: [0,t\wedge\tau^x_r]\times \Omega\longrightarrow \IR
\end{align*}
is a continuous local martingale \cite{driver}. Using (\ref{qbit1}), (\ref{qbit2}), the assumptions on $\ell$, the Burkholder-Davis-Gundy inequality (the latter to estimate $U^{(\ell)}$) and that $X^x$ takes values in a compact set on $[0,t\wedge\tau^x_r]$, the process $Y$ is in fact a true martingale by Lemma \ref{paaq}, in particular, $ Y$ has a constant expectation. Evaluating $Y_s$ at the times $s = 0$ and $s = t\wedge \tau^x_r$ and taking expectations, we get $\bE [Y_0] = \bE [ Y_{t\wedge \tau^x_r} ]$ so that
\begin{align*}
(    \nabla \mathrm{e}^{-\frac{t}{2}\vec{\Delta}_j  } \alpha(x),\xi) = -\bE \left[\left( Q^x_j(t\wedge \tau^x_r) \transport_{t\wedge \tau^x_r}^{x,-1} \mathrm{e}^{-\frac{(t-t\wedge \tau^x_r)}{2}}\alpha(X_{t\wedge \tau^x_r}^x) , U^{(\ell)}_{t\wedge \tau^x_r}\right)\right],
\end{align*}
which is the local Bismut derivative formula.

\end{proof}

\begin{Lemma}\label{ebni} For all $t>0$, $r>0$, $x\in M$, $\xi\in T^*_{x}M\otimes \Lambda^jT^*_xM$ there exists a process $\ell\in \IAA(x,r,t,\xi)$ such that for all $1\leq q<\infty$ and all constants $a\geq 0$ with $\mathrm{Ric}\geq-a$ in $B(x,r)$ one finds constants $C_{q,m},C_{a,q,m}<\infty$ satisfying 
\begin{align*}
|\ell|\leq |\xi|,\quad \mathbb{E}\left[\left(\int^{t\wedge \tau^x_r}_0 |\dot{\ell}_s|^2 \Id s\right)^{q/2}\right]^{1/q}\leq t^{-1/2}\mathrm{e}^{\frac{tC_{a,q,m}}{r}+\frac{tC_{q,m}}{r^2}}|\xi|. 
\end{align*}

\end{Lemma}

\begin{proof} It is well-known (cf. the proof of Corollary 5.1 in \cite{tw}) how to construct a bounded adapted process 
$$
k: [0,t] \times \Omega\longrightarrow \IR
$$
with paths in the Cameron-Martin space $W^{1,2}([0,t],\IR)$, such that   
$$
|k|\leq 1,\quad \mathbb{E}\left[\int_0^{t\wedge \tau^x_r }\abs{\dot k_s}^2 \Id s\right]<\infty,\quad  k_0 = 1,\quad k_s = 0\quad \text{for all $s \geq t\wedge\tau^x_r$},
$$
and
$$
\mathbb{E}\left[\left(\int^{t\wedge\tau^x_r}_0 |\dot{k}_s|^2 \Id s\right)^{q/2}\right]^{1/q}\leq t^{-1/2}\mathrm{e}^{\frac{tC_{a,q,m}}{r}+\frac{tC_{q,m}}{r^2}}.
$$
Thus we may simply set $\ell_s:=k_s \xi$.
 \end{proof}

\begin{Lemma}\label{katos} Assume $V_j\geq a$ for some constant $a\in\IR$, and let $\alpha\in \Gamma_{L^2}(M,\Lambda^j T^*M)$. Then for all $t> 0$ one has
$$
|\mathrm{e}^{-t \vec{\Delta}_j}\alpha|\leq \mathrm{e}^{-at}\mathrm{e}^{-t \Delta}|\alpha|.
$$
\end{Lemma}

\begin{proof} As already noted in the introduction, this semigroup domination is a well-known fact \cite{Ber}. Much more general statements, which do not require constant lower bounds on the potential, can be found in \cite{batu2} and are referred to as Kato-Simon inequality there.
\end{proof}

The following covariant Feynman-Kac formula is well-known in much more general situations \cite{batu,driver} to hold \emph{a.e. in $M$}; the point of the proof below (which is the usual one for compact $M$'s) is that it identifies the smooth representative of $\mathrm{e}^{-\frac{t}{2} \vec{\Delta}_j} \alpha$ \emph{pointwise on $M$}:

\begin{Lemma}[Covariant Feynman-Kac formula]\label{abedy} Assume $M$ is stochastically complete with $V_j\geq a$ for some constant $a\in\IR$. Then for all $t\geq 0$, $\alpha\in\Gamma_{L^2\cap L^{\infty}\cap C^{\infty}}(M,\Lambda^jT^*M)$, $x\in M$ one has
\begin{align}\label{posss}
 \mathrm{e}^{-\frac{t}{2}\vec{\Delta}_j  } \alpha(x)=\mathbb{E}\left[Q_{j}^x(t)  \transport_{t}^{x,-1} \alpha(X^x_t)\right].
\end{align}
\end{Lemma}

\begin{proof} Note that the statement of above formula includes that the right-hand side coincides for \emph{all} $x\in M$ and not only for $\mu$-a.e. $x\in M$ with the smooth representative of $\mathrm{e}^{-\frac{s}{2}\vec{\Delta}_j  }\alpha$. To prove the formula, we can assume that $t>0$. Then the process
\[
Y:[0,t]\times\Omega\longrightarrow \Lambda^jT^*_x M,\>\>Y_s:=Q^x_j(s)\transport^{x,-1}_s\mathrm{e}^{-\frac{t-s}{2}\vec{\Delta}_j}\alpha(X_s^x) \label{gelb}
\]
is a continuous local martingale. Under the stated assumptions, using Lemma \ref{katos} and $|Q^x_j(s)|\leq \mathrm{e}^{-as}$ $\mathbb{P}$-a.s. (by Gronwall\rq{}s lemma), one finds
\begin{align*}
 \left|Y_s\right|\leq \mathrm{e}^{2|a|t} \left\|\alpha\right\|_{\infty} \int \mathrm{e}^{-\frac{t-s}{2}\Delta}(X_s^x,y) \Id\mu (y)\leq  \mathrm{e}^{2|a|t}\left\|\alpha\right\|_{\infty}, 
\end{align*}
so that $Y$ is in fact a martingale by Lemma \ref{paaq}. Evaluating $Y_s$ at $s=0$ and $s=t$ and taking expectations proves the claim.
\end{proof}

\begin{Lemma}\label{apriori} Assume (\ref{wlp}). Then there exists a constant $C=C(A,m)>0$, such that for all $j\in\{1,\dots,m\}$, $t>0$, $x\in M$, $\alpha\in\Gamma_{L^2\cap C^{\infty}\cap L^{\infty}}$ one has 
$$
|\nabla \mathrm{e}^{-t\vec{\Delta}_j  } \alpha(x)|\leq C t^{-1/2}\mathrm{e}^{Ct}\left\|\alpha\right\|_{\infty}.
$$
\end{Lemma}

\begin{proof} In the sequel, $C(a,\dots)$ will denote a constant that only depends on $a,\dots,$ and which may differ from line to line. Let $t>0$, $r>0$, $x\in M$, $\xi\in T_{x}^*M\otimes \Lambda^jT^*_xM$ be arbitrary and pick $\ell\in \IAA(x,r,t,\xi)$ as in Lemma \ref{ebni}. We set
\begin{align*}
&\ell^{(1)}:= \int_0^{\bullet} Q^x_j(s)^{\dagger,-1} \Id\underline{X}^x_s \underline{Q}^x_j(s)^{\dagger}  \dot \ell_s,\\
&\ell^{(2)}:=\frac{1}{2} \int_0^{\bullet} Q^x_j(s)^{\dagger,-1}   \big(\transport_s^{x,-1}  \rho_j(X_s^x) \transport_s^x\big)^{\dagger} \underline{Q}^x_j(s)^{\dagger}  \ell_s \Id s.
\end{align*}
It follows from the covariant Feynman-Kac formula, the fact that the (inverse) damped parallel transport $Q_j  \transport^{-1}$ is a multiplicative functional (cf. equation (61) in \cite{batu}) and the strong Markov property of Brownian motion, that
$$
Q^x_j(t\wedge\tau^x_r)  \transport_{t\wedge \tau^x_r }^{x,-1}  \mathrm{e}^{\frac{-(t-t\wedge \tau^x_r)}{2} \vec{\Delta}_j}\alpha(X^x_{t\wedge \tau^x_r } ) = \mathbb{E}^{\mathcal{F}_{t\wedge \tau^x_r}}\left[  Q^x_j(t)  \transport_{t}^{x,-1}\alpha(X^x_t)\right],
$$
and so since $\ell^{(1)}_{t\wedge \tau^x_r }+\ell^{(2)}_{t\wedge \tau^x_r }$ is $\mathcal{F}_{t\wedge \tau^x_r}$-measurable, the law of total expectation gives 
\begin{align*}
&\mathbb{E}\left[\left(Q^x_j(t\wedge\tau^x_r)  \transport_{t\wedge \tau^x_r }^{x,-1}  \mathrm{e}^{\frac{-(t-t\wedge \tau^x_r)}{2} \vec{\Delta}_j}\alpha(X^x_{t\wedge \tau^x_r } ), \ell^{(1)}_{t\wedge \tau^x_r }+\ell^{(2)}_{t\wedge \tau^x_r }\right)\right]\\
&= \mathbb{E}\left[\mathbb{E}^{\mathcal{F}_{t\wedge \tau^x_r}}\left[\left(  Q^x_j(t)  \transport_{t}^{x,-1}\alpha(X^x_t),\ell^{(1)}_{t\wedge \tau^x_r }+\ell^{(2)}_{t\wedge \tau^x_r }\right)\right]\right]\\
&=\mathbb{E}\left[\left(  Q^x_j(t)  \transport_{t}^{x,-1}\alpha(X^x_t),\ell^{(1)}_{t\wedge \tau^x_r }+\ell^{(2)}_{t\wedge \tau^x_r }\right)\right].
\end{align*}

Thus (\ref{eq::nablaestimg}) implies
$$
(    \nabla \mathrm{e}^{-\frac{t}{2}\vec{\Delta}_j  } \alpha(x),\xi)=-\bE \left[\left( Q_{j}^x(t)  \transport_{t}^{x,-1}  \alpha(X_{t}^x) , \ell^{(1)}_{t\wedge \tau^x_r }+\ell^{(2)}_{t\wedge \tau^x_r }\right)\right],
$$
and we have
\begin{align*}
\mathbb{E}\left[|\ell^{(1)}_{t\wedge \tau^x_r }|\right]&\leq C \mathbb{E}\left[\left(\int^{t\wedge \tau^x_r }_0 |Q_s^{x,\mathrm{tr},-1} |^2| \underline{Q}_s^{x,\mathrm{tr}}|^2 |\dot{\ell}_s|^2 \Id s\right)^{1/2}\right]\\
&\leq C\mathrm{e}^{C(m)t}t^{-1/2}\mathrm{e}^{\frac{tC(A,m)}{r}+\frac{tC(m)}{r^2}}|\xi|,
\end{align*}
having used the Burkholder-Davis-Gundy inequality as well as (\ref{qbit1}) and (\ref{qbit2}). Moreover,  
$$
\mathbb{E}\left[|\ell^{(2)}_{t\wedge \tau^x_r }|\right]\leq    \mathrm{e}^{C(m)t}C(m,A)  |\xi|. 
$$
which follows from (\ref{qbit1}), (\ref{qbit2}), $|\ell|\leq |\xi|$, $|\rho|\leq C(m,A)$. Using once more (\ref{qbit1}), we can now estimate as follows
\begin{align*}
&|(    \nabla \mathrm{e}^{-\frac{t}{2} \vec{\Delta}_j  } \alpha(x),\xi)|\\
&\leq \bE \left[| Q_{j}^x(t)  | |  \alpha(X_{t}^x)| | \ell^{(1)}_{t\wedge \tau^x_r }|\right]+ \bE \left[| Q_{j}^x(t)|| \alpha(X_{t}^x)|| \ell^{(2)}_{t\wedge \tau^x_r }|\right]\\
&\leq |\xi|C(A,m)\mathrm{e}^{C(m,A)t}\left\|\alpha\right\|_{\infty}\Big(t^{-1/2}\mathrm{e}^{\frac{tC(A,m)}{r}+\frac{tC(m)}{r^2}}+1 \Big).
\end{align*}
Taking $r\to \infty$, we have managed to construct $C(A,m)<\infty$, such that for all $x\in M$, $t>0$, one has 
\begin{align*}
|    \nabla \mathrm{e}^{-\frac{t}{2} \vec{\Delta}_j  } \alpha(x)|\leq C(A,m) \mathrm{e}^{tC(A,m)}t^{-1/2}  \left\|\alpha\right\|_{\infty}.  
\end{align*}
\end{proof}

Being equipped with the latter a priori $L^{\infty}$ bound, we can now prove the global Bismut derivative formula. To this end, for fixed $x\in M$, $\zeta\in T^*_xM\otimes \Lambda^j_x T^*M$, $t>0$, define the continuous semimartingale 
\small{\begin{align*}
   & U(x,t,\zeta):[0,t]\times \Omega\longrightarrow T_{x}M ,\\
		& U_{\bullet}(x,t,\zeta):=  -\frac{1}{t}\int_0^{\bullet} Q^x_j(s)^{\dagger,-1} \Id\underline{X}^x_s \underline{Q}^x_j(s)^{\dagger}\zeta+\frac{1}{2t} \int_0^{\bullet} Q^x_j(s)^{\dagger,-1}   \big(\transport_s^{x,-1}  \rho_j(X_s^x) \transport_s^x\big)^{\dagger} \underline{Q}^x_j(s)^{\dagger}  (t-s)\zeta \Id s.
\end{align*}
}

\begin{Theorem}[Global Bismut derivative formula] Assume (\ref{wlp}). For every $t > 0$,  $x\in M$, $\xi\in T_{x}^*M\otimes \Lambda^jT_{x}M$, and every $\alpha\in \Gamma_{L^2\cap C^{\infty}\cap L^{\infty}}(M,\Lambda^jT^*M)$ one has
\begin{align}\label{eq::nablaestimg2}
 (  \nabla  \mathrm{e}^{-\frac{t}{2}\vec{\Delta}_j } \alpha(x)  , \xi   ) = -\bE  \left[ \big( Q^x_j(t)  \transport_{t}^{x,-1}  \alpha(X^x_{t} )  , U_t(x,t,\xi) \big)\right].
\end{align}
 \end{Theorem}
 
\begin{proof} For $0\leq s\leq t$, let $\ell(s):=(1-\frac{s}{t})\zeta$. As before, the process

\begin{align*} 
   Y &:= \big( \tilde Q^x_j   \transport ^{x,-1}  \nabla \mathrm{e}^{-\frac{t-\bullet}{2}\vec{\Delta}_j}  \alpha(X^x ) , \ell\big) - \big( Q^x_j   \transport^{x,-1} \mathrm{e}^{-\frac{t-\bullet}{2}\vec{\Delta}_j} \alpha(X^x ) , U(x,t,\zeta) \big) \\
   &: [0,t]\times \Omega\longrightarrow \IR
\end{align*}
is a continuous local martingale (cf \cite{driver}). Using Lemma \ref{apriori} we have 
$$
|Y_s|\leq C(A,m) \mathrm{e}^{C(A,m)t}|\zeta|\left\|\alpha\right\|_{\infty}\big(\frac{(t-s)^{1/2}}{t}+|U_s|\big),
$$
and using the Burkholder-Davis-Gundy inequality as well as (\ref{qbit1}) and (\ref{qbit2}), one easily finds
$$
\mathbb{E}\left[\sup_{s\in [0,t]}|Y_s|\right]<\infty,
$$
showing that $Y$ is a martingale, and the global Bismut derivative formula follows from $\mathbb{E}\left[Y_0\right]=\mathbb{E}\left[Y_t\right]$.
\end{proof}

\section{Proof of Theorem \ref{main}}\label{beweis}
With $\mathbb{E}^{x,y}_t\left[\bullet\right]$ denoting integration with respect to the Brownian bridge measures \cite{bbb}, the global Bismut derivative formula together with the disintegration property
$$
\mathbb{E}[\Psi(X^x)]= \int_M\mathrm{e}^{-\frac{t}{2}\Delta }(x,y) \mathbb{E}^{x,y}_t[\Psi(X^x)]\Id\mu(y),
$$
valid for all Borel-measurable (vector-valued) functions $\Psi$ on the space of continuous paths $C([0,t], M)$, one has
\begin{align*}
&\left|\nabla \mathrm{e}^{-\frac{t}{2}\vec{\Delta}_j}(x,y)\right|\\
&\leq \frac{1}{t}\mathrm{e}^{-\frac{t}{2}\Delta }(x,y)\mathbb{E}^{x,y}_t\left[|Q^x_j(t)| \left|\int_0^{t}   Q_j^x(s)^{\dagger,-1} \Id \underline{X}^x_s \underline{Q}^x_j(x)^{\dagger}  \right|\right]\\
& \quad +\mathrm{e}^{-\frac{t}{2}\Delta }(x,y) \frac{1}{2t}\mathbb{E}^{x,y}_t\left[|Q^x_j(t)| \int_0^{t}  |  Q^x_j(s)^{\dagger,-1} | | \rho(X_s^x) | (t-s) \Id s\right]\\
&\leq \frac{1}{t}\mathrm{e}^{-\frac{t}{2}\Delta }(x,y) \mathrm{e}^{C(A,m)t}\mathbb{E}^{x,y}_t\left[  \left|\int_0^{t}   Q^x_j(s)^{\dagger,-1} \Id X^x_s \underline{Q}_j^x(s)^{\dagger} \right|\right]+ C(A,m)\mathrm{e}^{C(A,m)t}\mathrm{e}^{-\frac{t}{2}\Delta }(x,y). 
\end{align*}
Furthermore, using the time reversal property and the defining relation of the Brownian bridge \cite{bbb} we have
\begin{align*}
&\mathrm{e}^{-\frac{t}{2}\Delta }(x,y) \mathbb{E}^{x,y}_t\left[  \left|\int_0^{t}   Q^x_j(s)^{\dagger,-1} \Id X^x_s \underline{Q}_j^x(s)^{\dagger}  \right|\right]\\
&= \mathbb{E}\left[ \mathrm{e}^{-\frac{t}{2}\Delta }(X^x_{t/2},y)  \left|\int_0^{t/2}   Q^x_j(s)^{\dagger,-1} \Id X^x_s \underline{Q}_j^x(s)^{\dagger}  \right|\right]\\
&+\mathbb{E}\left[ \mathrm{e}^{-\frac{t}{2}\Delta }(X^y_{t/2},x)  \left|\int_0^{t/2}   Q^x_j(s)^{\dagger,-1} \Id X^x_s \underline{Q}_j^x(s)^{\dagger}  \right|\right],
\end{align*}
and using Cauchy-Schwarz and estimating the stochastic integral using Burkholder-Davis-Gundy, keeping in mind that $X^x|_{[0,t]}$ is a semimartingale under the Brownian bridge measure \cite{batu},
\begin{align*}
&\mathbb{E} \left[ \mathrm{e}^{-\frac{t}{2}\Delta }(X^x_{t/2},y) \left|\int_0^{t/2}   Q_j^x(s)^{\dagger,-1} \Id X^x_s \underline{Q}_j^x(s)^{\dagger}  \right|\right]\\
&\leq \mathbb{E} \left[ \mathrm{e}^{-\frac{t}{2}\Delta }(X^x_{t/2},y)^2  \right]^{1/2}\mathbb{E} \left[\left|\int_0^{t/2}   Q_j^x(s)^{\dagger,-1} \Id X^x_s \underline{Q}_j^x(s){\dagger} \right|^2\right]^{1/2}\\
&\leq C(A,m)t^{1/2}\mathrm{e}^{C(A,m)t}\left(\int \mathrm{e}^{-\frac{t}{2}\Delta }(x,z)\mathrm{e}^{-\frac{t}{2}\Delta }(z,y) \mathrm{e}^{-\frac{t}{2}\Delta }(z,y) \Id\mu(z)\right)^{1/2}\\
&\leq C(A,m)\mu(B(y,\sqrt{t/2}))^{-1/2} t^{1/2}\mathrm{e}^{C(A,m)t}\mathrm{e}^{-t\Delta }(x,y)^{1/2}\\
&\leq C(A,m)\mu(B(y,\sqrt{t/2}))^{-1}t^{1/2}\mathrm{e}^{C(A,m)t}\mathrm{e}^{-C(A,m)\frac{\varrho(x,y)^2}{t}},
\end{align*}
where we have used the Li-Yau estimate
$$
\mathrm{e}^{-s\Delta }(x_1,x_2)\leq \mu(B(x_1,\sqrt{s}))^{-1}\mathrm{e}^{- C_{A,m} \frac{\varrho(x_1,x_2)^2}{s}}\mathrm{e}^{C(A,m)s}\quad\text{ for all $s>0$, $x_1,x_2\in M$},
$$
twice, and 
$$
\int \mathrm{e}^{-\frac{t}{2}\Delta }(x,z)\mathrm{e}^{-\frac{t}{2}\Delta }(z,y)   \Id\mu(z)= \mathrm{e}^{-t\Delta }(x,y).
$$
Likewise, we have
\begin{align*}
&\mathbb{E} \left[ \mathrm{e}^{-\frac{t}{2}\Delta }(X^x_{t/2},y) \left|\int_0^{t/2}   Q_j^x(s)^{\dagger,-1} \Id X^x_s \underline{Q}_j^x(s)^{\dagger}  \right|\right]\\
&\leq C(A,m)\mu(B(x,\sqrt{t/2}))^{-1}t^{1/2}\mathrm{e}^{C(A,m)t}\mathrm{e}^{-C(A,m)\frac{\varrho(x,y)^2}{t}},
\end{align*}
so that with local doubling \eqref{eq:LVD} we arrive at the desired estimate.

\section{Proof of Corollary \ref{cor:gradient}}\label{beweis2}

\subsection{Proof of Corollary \ref{cor:gradient}.I)}

Let $K_t$ be the integral operator with Gaussian integral  kernel 
$$
k_t(x,y):=\frac{1}{\mu(B(x,\sqrt{t}))}\mathrm{e}^{-\frac{-c\varrho^2(x,y)}{t}}
$$
according to the estimate of Theorem \ref{main}, it is enough to prove that there is a constant $C>0$ such that for all $p\in [1,+\infty]$,
\begin{equation}\label{eq:Lpaction}
||K_t||_{p,p}\leq C\mathrm{e}^{Ct}
\end{equation}
By interpolation, it is enough to prove \eqref{eq:Lpaction} for $p=1$ and $p=\infty$. However,
$$||K_t||_{1,1}=\sup_{y\in M}\int_M k_t(x,y)\, \Id\mu(x),$$
and likewise,
$$||K_t||_{\infty,\infty}=\sup_{x\in M}\int_M k_t(x,y)\, \Id\mu(y).$$
The volume comparison inequality \eqref{eq:VC} with small enough $\epsilon$ implies that there exist constants $C_1,C_2>0$ such that for all $x,y\in M$ and all $t>0$,
$$\frac{1}{C_2\mu(B(y,\sqrt{t}))}\mathrm{e}^{C_2t}\mathrm{e}^{-\frac{-\varrho^2(x,y)}{C_2t}}\leq \frac{1}{\mu(B(x,\sqrt{t}))}\mathrm{e}^{-\frac{-\varrho^2(x,y)}{Ct}}\leq \frac{C_1}{\mu(B(y,\sqrt{t}))}\mathrm{e}^{C_1t}\mathrm{e}^{-\frac{-\varrho^2(x,y)}{C_1t}}$$
As a consequence, $||K_t||_{1,1}\leq C\mathrm{e}^{Ct}$ follows from $||K_t||_{\infty,\infty}\leq C\mathrm{e}^{Ct}$. Hence, \eqref{eq:Lpaction} will follow from the estimate: there is a constant $C>0$ such that

\begin{equation}\label{eq:Linftyaction}
\sup_{x\in M}\int_M k_t(x,y)\, \Id\mu(y)\leq C\mathrm{e}^{Ct}
\end{equation}
Let $x\in M$, $A_0=2B(x,\sqrt{t})$ and, for all $i\geq 1$, $A_i=2^{i+1}B(x,\sqrt{t})\setminus 2^iB(x,\sqrt{t})$. Then, using \eqref{eq:LVD}, one gets

\begin{eqnarray*}
\int_M k_t(x,y)\, \Id\mu(y) & = & \sum_{i=0}^\infty \int_{A_i}k_t(x,y)\, \Id\mu(y)\\
&\leq & C\sum_{i=0}^\infty \frac{\mu(A_i)}{\mu(B(x,\sqrt{t}))}\mathrm{e}^{-c4^i}\\
&\leq & C\sum_{i=0}^\infty \frac{\mu(B(x,2^{i+1}\sqrt{t}))}{\mu(B(x,\sqrt{t}))}\mathrm{e}^{-c4^i}\\
&\leq & C\sum_{i=0}^\infty 2^{im} \mathrm{e}^{C2^{i+1}\sqrt{t}}\mathrm{e}^{-c4^i}\\
\end{eqnarray*}
Using the elementary inequality $\mathrm{e}^{C2^{i+1}\sqrt{t}}=\mathrm{e}^{\frac{C}{\epsilon}t}\mathrm{e}^{-\epsilon 4^i}$ with $\epsilon=\frac{c}{2}$, we arrive to

\begin{eqnarray*}
\int_M k_t(x,y)\, \Id\mu(y) &\leq & C\sum_{i=0}^\infty 2^{im}\mathrm{e}^{Ct}\mathrm{e}^{-C4^i}\\
&\leq & C\mathrm{e}^{Ct},
\end{eqnarray*}
which completes the proof of \eqref{eq:Linftyaction}.

\subsection{Proof of Corollary \ref{cor:gradient}.II)}

By Theorem \ref{main} and \eqref{eq:VC}, given $\gamma>0$, we get
\begin{align*}
&\int |\nabla \mathrm{e}^{-t \vec{\Delta}_j}(x,y)|^p \mathrm{e}^{\frac{\gamma \varrho(x,y)^2}{t}}\Id\mu(x)\\
&\leq C(A,m,p)\mathrm{e}^{C(A,m)t} t^{-p/2}\mu(B(y,\sqrt{t}))^{-p}\int \mathrm{e}^{(\gamma-C_{m,p,A})\frac{\varrho(x,y)^2}{t}}\Id\mu(x)\\
&\leq C(A,m,p)\mathrm{e}^{C(A,m)t}t^{-p/2}\mu(B(y,\sqrt{t}))^{-p}\sum^{\infty}_{i=1}\int_{B(y,i\sqrt{t})\setminus B(y,(i-1)\sqrt{t})} \mathrm{e}^{(\gamma-C_{m,p,A})\frac{\varrho(x,y)^2}{t}}\Id\mu(x)\\
&\quad +C(A,m,p)\mathrm{e}^{C(A,m)t}t^{-p/2} \mu(B(y,\sqrt{t}))^{-p-1},
\end{align*}
where we have chosen $\gamma<C(A,m,p)$. Finally, using local doubling \eqref{eq:LVD}, letting $\gamma':= C(A,m,p)-\gamma>0$, we have
\begin{align*}
&\sum^{\infty}_{i=1}\int_{B(y,i\sqrt{t})\setminus B(y,(i-1)\sqrt{t})} \mathrm{e}^{-\gamma'\frac{\varrho(x,y)^2}{t}}\Id\mu(x)\\
&\leq \mu(B(y,\sqrt{t}))\sum^{\infty}_{i=1} \frac{\mu(B(y,i\sqrt{t}))}{\mu(B(y,\sqrt{t}))} \mathrm{e}^{-\gamma'(i-1)^2}\\
&\leq \mu(B(y,\sqrt{t}))\sum^{\infty}_{i=1} j^m \mathrm{e}^{-\gamma'(i-1)^2+C(A,m)i}<\infty,
\end{align*}
completing the proof.

\section{Proof of Corollary \ref{cor:Riesz} and Theorem \ref{thm:riesz-bakry}}\label{RE}

In this section, we explain how one can use the heat kernel estimates \eqref{eq:vUE}, \eqref{eq:CUE}, \eqref{eq:dj} and \eqref{eq:d*j} in order to get results for the Riesz transforms $\nabla (\vec{\Delta}_j+\lambda)^{-1/2}$ and $(\Id_j+\Id^\dagger_{j-1})(\vec{\Delta}_j+\lambda)^{-1/2}$ (i.e. prove Corollary \ref{cor:Riesz} and Theorem \ref{thm:riesz-bakry} respectively). The idea of the proof is to follow closely the proof of \cite[Theorem 1.2]{cd}, where a result for the localized scalar Riesz transform $\Id(\Delta+\lambda)^{-1/2}$ is proved. The proof is based on the Calder\'on-Zygmund decomposition and kernel estimates, which follow from the assumed heat kernel estimates \eqref{eq:vUE}, \eqref{eq:CUE}, \eqref{eq:dj} and \eqref{eq:d*j}. However we feel that in the proof of \cite[Theorem 1.2]{cd} the issue of localization may have been overlooked a little: there, it is wrongly asserted that \eqref{eq:LVD} implies that every open ball of radius $1$ in $M$ is a doubling space, with a doubling constant that can be chosen independently of the ball; actually, this property depends on the geometry of balls, and not only on the validity of \eqref{eq:LVD} in the whole $M$, and we don't see why it should hold in the context of \cite[Theorem 1.2]{cd}. In order to clarify the matter, we decided to give full proofs for the localization procedure that we use. The first ingredient needed in our proof is a localized Calder\'on-Zygmund decomposition $f=g+b$ for a smooth section $f\in \Gamma(M,\Lambda^j T^*M)$ which has support inside a ball $B=B(x,1)$. This decomposition holds thanks to the local doubling assumption \eqref{eq:LVD}. More precisely, the version of the Calder\'on-Zygmund decomposition we need is the following: 

\begin{Lemma}\label{lem:CZ}

Let $\IEE\to M$ be a Riemannian vector bundle, where $M$ is locally doubling. Then there is a constant $C>0$, which depends only on the local doubling constant, with the following property: for every ball $B=B(x,1)$, every $u\in \Gamma_{C^\infty}(M,\IEE)$ with support inside $B$, and every $0<\lambda<\frac{C}{\mu(B)}\int_B |u|$, there exists a countable collection of balls $(B_i)_{i\in I}$, of integrable sections $(b_i)_{i\in I}$ in $\Gamma_{L^1}(M,\IEE)$ and a section $g\in \Gamma_{L^\infty}(M,\IEE)$ such that:

\begin{enumerate}

\item $u=g+\sum_{i\in I} b_i$ a.e.

\item the balls $(B_i)_{i\in I}$ have the finite intersection property: there is $N\in \N$ such that for every $i\in\N$,

$$\mathrm{Card}\{j\in \N\,:\,B_i\cap B_j\neq \emptyset\}\leq N.$$

\item $\sum_{i\in I}\mu(B_i)\leq \frac{C}{\lambda}\int_B |u|$.

\item $|g|\leq \lambda$ a.e.

\item For all $i\in I$, $b_i$ has support inside $B_i$, and $\int_{B_i}|b_i|\leq C\lambda \mu(B_i)$.

\end{enumerate}

\end{Lemma}
Furthermore, as a consequence of (2), (3) and (5), there holds for some constant $C$:

\begin{equation}\label{eq:g}
||g||_1\leq C||u||_1
\end{equation}
The proof of this version of the Calder\'on-Zygmund decomposition closely follows the classical one, with three differences: firstly, since one has only local doubling but not doubling, one has to use a modified maximal function $\mathfrak{M}$, defined as follows:

$$\mathfrak{M}u(x):=\sup_{B\ni x\,:\,r(B)\leq 8}\frac{1}{\mu(B)}\int_B|u|,$$
where $r(B)$ denotes the radius of the ball $B$. The particular value $8$ in the definition of $\mathfrak{M}$ is chosen for later technical purposes (see the proof of Lemma \ref{lem:max}). Note that local doubling implies that $\mathfrak{M}$ is weak type $(1,1)$ and bounded on $L^p$ for all $p\in (1,\infty]$, as follows from a careful inspection of the proof of \cite[Theorem 1 p. 13]{St} and the fact that the definition of $\mathfrak{M}$ involves only balls with bounded radii. Secondly, in the Calder\'on-Zygmund decomposition localized in the ball $B$, the balls $B_i$ do not have to be included inside the ball $B$, only inside $2B$. Lastly, the fact that we are dealing here with sections of a vector bundle instead of mere functions: this does not create any real difficulty and the standard arguments apply {\em mutatis mutandis} if one puts norms instead of absolute values everywhere it is needed. A detailed proof of Lemma \ref{lem:CZ} is presented in Appendix C.\\

Let us now present the main steps of the proof of  Corollary \ref{cor:Riesz} and Theorem \ref{thm:riesz-bakry}, following closely the approach of \cite[Theorem 1.2]{cd}. Let $T$ be either $\nabla (\vec{\Delta}_j+\lambda)^{-1/2}$ or $(\Id_j+\Id^\dagger_{j-1})(\vec{\Delta}_j+\lambda)^{-1/2}$. We start with boundedness of $T$ on $L^2$:

\begin{Lemma}\label{lem:L2} For all $\kappa>0$ the operator $(\Id_j+\Id^\dagger_{j-1})(\vec{\Delta}_j+\kappa)^{-1/2}$, originally defined on $\Gamma_{C_c^\infty}(M,\Lambda^jT^*M)$, extends to a bounded operator on $\Gamma_{L^2}(M,\Lambda^jT^*M)$ with 
$$
\left\|(\Id_j+\Id^\dagger_{j-1})(\vec{\Delta}_j+\kappa)^{-1/2}\right\|_{2,2}\leq 1.
$$
If $||\mathrm{Riem}||_\infty\leq A<\infty$, then for all $\kappa>0$ the operator $\nabla(\vec{\Delta}_j+\kappa)^{-1/2}$, originally defined on $\Gamma_{C_c^\infty}(M,\Lambda^jT^*M)$, extends to a bounded operator on $\Gamma_{L^2}(M,\Lambda^jT^*M)$ with 
$$
\left\|\nabla(\vec{\Delta}_j+\kappa)^{-1/2}\right\|_{2,2}\leq C,
$$
where $C$ only depends on $\lambda$, $A$, $m$.

\end{Lemma}

\begin{proof} Since $\Gamma_{C_c^\infty}(M,\Lambda^jT^*M)$ is dense in $\Gamma_{L^2}(M,\Lambda^jT^*M)$, it is enough to show for any $f\in \Gamma_{C_c^\infty}(M,\Lambda^jT^*M)$, 
\begin{equation}\label{eq:RL2}
||(\Id_j+\Id^\dagger_{j-1})(\vec{\Delta}_j+\kappa)^{-1/2}f||_2\leq ||f||_2,
\end{equation}
and
\begin{equation}\label{eq:RCL2}
||\nabla (\vec{\Delta}_j+\kappa)^{-1/2}f||_2\leq C||f||_2.
\end{equation}

The first estimate is a simple consequence of the functional calculus: since the Dirac operator $D:=\Id+\Id^\dagger$ acting on smooth, compactly supported differential forms, is essentially self-adjoint on $M$, it follows that for all $g$ in the domain $D^2=\vec{\Delta}$ (which is included in the domain of $D$),
\begin{eqnarray*}
||(\Id+\Id^\dagger)g||_2^2 &=& \langle Dg,Dg \rangle\\
&=& \langle D^2g,g\rangle\\
&\leq &   \langle (D^2+\kappa)g,g\rangle \\
&\leq & ||(D^2+\kappa)^{1/2}g||_2^2.
\end{eqnarray*}
Apply the above inequality to $g=(D^2+\kappa)^{-1/2}f$, which is the domain of $D^2$ by functional calculus, we obtain \eqref{eq:RL2} with $C=1$.
Let us now prove \eqref{eq:RCL2}. Recall that since $M$ is complete, the operator $\nabla^\dagger \nabla$ acting on smooth compactly supported differential forms is essentially self-adjoint, associated with the quadratic form $(u,v)\mapsto \langle \nabla u,\nabla v\rangle$. In particular, if $g\in \Gamma_{L^2}(M,\Lambda^jT^*M)$ is in the domain of $\nabla^\dagger\nabla$, then
$$\langle \nabla g,\nabla g \rangle = \langle \nabla^\dagger\nabla g,g\rangle.$$
Hence, for such a $g$, using that $||V_j||_\infty\leq A'$, where $A'=A'(A,m)<\infty$,
\begin{eqnarray*}
||\nabla g||_2^2 &=& \langle \nabla g,\nabla g \rangle\\
&=& \langle \nabla^\dagger\nabla g,g\rangle\\
&\leq &   \langle (\nabla^\dagger\nabla+V_j+\kappa)g,g\rangle+A'||g||_2 \\
&\leq & ||(\vec{\Delta}_j+\kappa)^{1/2}g||_2^2+A'||g||_2^2.
\end{eqnarray*}
Take $g=(\vec{\Delta}_j+\kappa)^{-1/2}f$, which is in the domain of $\vec{\Delta}_j$: indeed, writing 
$$
f=(\vec{\Delta}_j+1)^{-1}(\vec{\Delta}_j+1)f,
$$
which can be done, since being smooth and compactly supported, $f$ is in the domain of $\vec{\Delta}_j$, one has
$$
g=(\vec{\Delta}_j+1)^{-1}(\vec{\Delta}_j+1)^{-1/2}(\vec{\Delta}_j+1)f,
$$
so that $g$ is in the domain of $\vec{\Delta}_j$ by functional calculus. It follows that
\begin{eqnarray*}
||\nabla  (\vec{\Delta}_j+\kappa)^{-1/2}f||_2^2 &\leq & ||f||_2^2 +A'||(\vec{\Delta}_j+\kappa)^{-1/2}f||_2^2\\
&\leq & (\frac{A'}{\kappa}+1)||f||_2^2,
\end{eqnarray*}
where we have used that $||(\vec{\Delta}_j+\kappa)^{-1/2}||_{2,2}\leq \kappa^{-1/2}$ by functional calculus. This proves \eqref{eq:RCL2}.

\end{proof}

Let us now come to the actual proof of Corollary \ref{cor:Riesz} and Theorem \ref{thm:riesz-bakry}: given the result of Lemma \ref{lem:L2} and using interpolation, one sees that it is enough to prove that $T$ is bounded from $\Gamma_{L^1}(M,\Lambda^j T^*M)$ to the space of weakly integrable sections $\Gamma_{L^1_w}(M,\Lambda^j T^*M)$, that is: one can find a constant $C>0$ such that for all $f\in \Gamma_{L^1}(M,\Lambda^jT^*M)$ and all $\lambda>0$,
\begin{equation}\label{eq:L1w}
\mu(\{x\in M\,:\,|Tf|(x)>\lambda \})\leq \frac{C}{\lambda}||f||_1.
\end{equation}
By a density argument, it is enough to prove it for $f$ smooth with compact support. So, take such an $f$, and fix $\lambda>0$. Take $(x_j)_{j\in\mathbb{N}}$ a maximal $1$-separated subset, hence the balls $B(x_j,1)$ cover $M$, while the balls $B(x_j,\frac{1}{2})$ are disjoint. Local doubling then implies that the balls $B(x_j,1)$ have the finite intersection property. Let $(\phi_j)_{j\in \N}$ be a smooth partition of unity associated to the covering of $M$ by the balls $B(x_j,1)$, and let $f_j:=\phi_jf$. The fact that the covering has the finite intersection property implies that for some constant $C>0$,
$$C^{-1}||f||_1\leq \sum_{j\in\N}||f_j||_1\leq C||f||_1.$$
Hence, it is enough to prove \eqref{eq:L1w} for $f_j$ (with a constant independent of $j$ and $f$). In what follows, we therefore assume that $j\in \N$ is fixed, and let $u=f_j$ and $B=B(x_j,1)$. We have two cases, according to whether 
$$\lambda\leq \frac{C}{\mu(B)}\int_B |u|$$
or not (here, $C$ is the constant in Lemma \ref{lem:CZ}). We first treat the case where $\lambda\leq  \frac{C}{\mu(B)}\int_B |u|$, for which there are two steps: first, show that
\begin{equation}\label{eq:small1}
\mu(\{x\in 2B\,:\,|Tf|(x)>\lambda\})\leq \frac{C}{\lambda}||u||_1,
\end{equation}
and then show that
\begin{equation}\label{eq:small2}
\mu(\{x\in M\setminus 2B\,:\,|Tf|(x)>\lambda\})\leq \frac{C}{\lambda}||u||_1.
\end{equation}
For \eqref{eq:small1}, notice that $\{x\in 2B\,:\,|Tf|(x)>\lambda\}\subset 2B$, therefore 
\begin{eqnarray*}
\mu(\{x\in 2B\,:\,|Tf|(x)>\lambda\}) &\leq & \mu(2B)\\
&\leq & C\mu(B)\\
&\leq & \frac{C}{\lambda} ||u||_1,
\end{eqnarray*}
where we have used successively \eqref{eq:LVD} and the assumption on $\lambda$. This proves \eqref{eq:small1}.

Now let us prove \eqref{eq:small2}. By the Markov inequality, we see that \eqref{eq:small2} follows from the $L^1$ estimate:

\begin{equation}\label{eq:small3}
\int_{M\setminus 2B}|Tu|(x)\, \Id\mu(x)\leq C||u||_1.
\end{equation}
In turns, \eqref{eq:small3} can be proved as in \cite[p. 1163]{cd}, using the heat kernel estimates \eqref{eq:CUE}, \eqref{eq:dj} and \eqref{eq:d*j} respectively.

\medskip

Now, we deal with the case $\lambda >\frac{C}{\mu(B)}\int_B |u|$. In this case one can use the Calder\'on-Zygmund decomposition $u=g+\sum_{i\in I}b_i$ from Lemma \ref{lem:CZ}. Let $r_i$ be the radius of $B_i$, and let $t_i=r_i^2$. Then, write

$$Tu=Tg+\sum_{i\in I}T\chi_{3B}\mathrm{e}^{-t_i\vec{\Delta}_j}b_i+\sum_{i\in I}T(1-\mathrm{e}^{-t_i\vec{\Delta}_j})b_i+ \sum_{i\in I}T\chi_{M\setminus 3B}\mathrm{e}^{-t_i\vec{\Delta}_j}b_i,$$
where we recall that $\chi_A$ denotes the indicator function of the Borel set $A$. The weak $L^1$ estimate \eqref{eq:L1w} will follow from the four estimates:
\begin{equation}\label{eq:Riesz1}
\mu(\{x\in M\,:\,|Tg|>\frac{\lambda}{4}\})\leq \frac{C}{\lambda}||f||_1,
\end{equation}

\begin{equation}\label{eq:Riesz2}
\mu(\{x\in M\,:\,|\sum_{i\in I}T\chi_{3B}\mathrm{e}^{-t_i\vec{\Delta}_j}b_i|>\frac{\lambda}{4}\})\leq \frac{C}{\lambda}||f||_1,
\end{equation}

\begin{equation}\label{eq:Riesz3}
\mu(\{x\in M\,:\,|  \sum_{i\in I}T(1-\mathrm{e}^{-t_i\vec{\Delta}_j})b_i  |>\frac{\lambda}{4}\})\leq \frac{C}{\lambda}||f||_1,
\end{equation}

\begin{equation}\label{eq:Riesz4}
\mu(\{x\in M\,:\,|\sum_{i\in I}T\chi_{M\setminus 3B} \mathrm{e}^{-t_i\vec{\Delta}_j}b_i  |>\frac{\lambda}{4}\})\leq \frac{C}{\lambda}||f||_1.
\end{equation}
Firstly, the fact that $|g|\leq C\lambda$ a.e. and that $T$ is bounded on $L^2$ leads to
\begin{eqnarray*}
\mu(\{x\in M\,:\,|Tg|>\frac{\lambda}{4}\}) &\leq & \frac{16}{\lambda^2} ||g||_2^2\\
&\leq & \frac{16}{\lambda^2}||g||_\infty ||g||_1\\
&\leq & \frac{C}{\lambda} ||u||_1,
\end{eqnarray*} 
which shows \eqref{eq:Riesz1}. Concerning \eqref{eq:Riesz2}, the same argument using the $L^2$ boundedness of $T$ shows that \eqref{eq:Riesz2} will follow from the $L^2$ estimate:
\begin{equation}\label{eq:Riesz21}
\left\|  \sum_{i\in I}\chi_{3B}\mathrm{e}^{-t_i\vec{\Delta}_j}b_i  \right\|_2^2 \leq C\lambda  ||u||_1
\end{equation}
Let $B_i=B(y_i,r_i)$. Point (5) of the Calder\'on-Zygmund decomposition together with the heat kernel estimate \eqref{eq:vUE} and the fact that $t_i\leq 2$ (since $B_i\subset 2B$) imply that
\begin{eqnarray*}
|\mathrm{e}^{-t_i\vec{\Delta}_j}b_i|(x) &\leq & C\lambda  \mu(B_i) \frac{\mathrm{e}^{-\frac{\varrho^2(x,y_i)}{Ct_i}}}{V(x,\sqrt{t_i})} \\
&\leq & C\lambda \int_{M}\frac{\mathrm{e}^{-\frac{\varrho^2(x,y)}{Ct_i}}}{V(x,\sqrt{t_i})}\chi_{B_i}(y)\, \Id\mu(y) \\
&\leq & C\lambda  \int_{M}\left(1+\frac{\varrho(x,y)}{\sqrt{t_i}}\right)^m\frac{\mathrm{e}^{-\frac{\varrho^2(x,y)}{Ct_i}}}{V(y,\sqrt{t_i})}\chi_{B_i}(y)\, \Id\mu(y) \\
&\leq &  C\lambda  \int_{M}\frac{\mathrm{e}^{-\frac{\varrho^2(x,y)}{Ct_i}}}{V(y,\sqrt{t_i})}\chi_{B_i}(y)\, \Id\mu(y),
\end{eqnarray*}
where in the one before last line we have used local doubling \eqref{eq:LVD} with $t_i\leq 2$, and where we have set
$$
V(z,r):=\mu(B(z,r)),\quad z\in M, r>0.
$$

In order to prove \eqref{eq:Riesz21}, it is then enough to prove that

\begin{equation}\label{eq:Riesz22}
\left\|  \chi_{3B}\sum_{i\in I} \int_{M}\frac{\mathrm{e}^{-\frac{\varrho^2(\cdot,y)}{Ct_i}}}{V(y,\sqrt{t_i})}\chi_{B_i}(y)\, \Id\mu(y) \right\|_2^2 \leq \frac{C}{\lambda}  ||u||_1.
\end{equation}
To estimate the above $L^2$ norm, we dualize against $v\in \Gamma_{L^2}(M,\Lambda^jT^*M)$ with support inside $3B$; we have by Fubini,

$$
\int_{M\times M} \sum_{i\in I} \frac{\mathrm{e}^{-\frac{\varrho^2(x,y)}{Ct_i}}}{V(y,\sqrt{t_i})}\chi_{B_i}(y) v(x)\, \Id\mu(x) \, \Id\mu(y) =\sum_{i\in I} \int_{B_i} \frac{1}{V(y,\sqrt{t_i})}\left(\int_{3B}  \mathrm{e}^{-\frac{\varrho^2(x,y)}{Ct_i}} v(x)\,\Id\mu(x)\right)\,\Id\mu(y).
$$

Next we proof that for every $i\in I$ and $y\in B_i$ one has
\begin{align}\label{lem:max}
\frac{1}{V(y,\sqrt{t_i})}\left(\int_{3B}  \mathrm{e}^{-\frac{\varrho^2(x,y)}{Ct_i}} v(x)\, \Id\mu(x)\right)\leq C\mathfrak{M}u(y).
\end{align}
Indeed, for $k\in \N$, let $A_0=B_i$ and
$$A_k=\{x\in3B\,:\,2^k\sqrt{t_i}\leq \varrho(x,y)\leq 2^{k+1}\sqrt{t_i}\},\quad k\geq 1.$$
Also, let $N\in \N$ be the smallest integer so that $2^{N+1}\sqrt{t_i}\geq 4$. Then, 

\begin{eqnarray*}
\frac{1}{V(y,\sqrt{t_i})}\left(\int_{3B}  \mathrm{e}^{-\frac{\varrho^2(x,y)}{Ct_i}} v(x)\, \Id\mu(x)\right) & = & \sum_{k=0}^\infty  \frac{1}{V(y,\sqrt{t_i})}\int_{A_k} \mathrm{e}^{-\frac{\varrho^2(x,y)}{Ct_i}} v(x)\, \Id\mu(x)\\
&\leq & \sum_{k=0}^N \frac{V(y,2^{k+1}\sqrt{t_i})}{V(y,\sqrt{t_i})}\mathrm{e}^{-c2^k} \frac{1}{V(y,2^{k+1}\sqrt{t_i})}\int_{B(y,2^{k+1}\sqrt{t_i})} |v|.
\end{eqnarray*}
By definition of $N$, we have for every $k\leq N$, $2^{k+1}\sqrt{t_i}\leq 8$, therefore by local doubling,

$$\frac{V(y,2^{k+1}\sqrt{t_i})}{V(y,\sqrt{t_i})}\leq C2^{km},$$
and it follows by definition of $\mathfrak{M}$ that

\begin{eqnarray*}
\frac{1}{V(y,\sqrt{t_i})}\left(\int_{3B}  \mathrm{e}^{-\frac{\varrho^2(x,y)}{Ct_i}} v(x)\, \Id\mu(x)\right) & \leq & \sum_{k=0}^N 2^{km}\mathrm{e}^{-c2^k} \mathfrak{M}v(y)\\
& \leq & \sum_{k=0}^\infty 2^{km}\mathrm{e}^{-c2^k} \mathfrak{M}v(y)\\
&\leq & C\mathfrak{M}v(y),
\end{eqnarray*}
and (\ref{lem:max}) is proved.

According to the remark made immediately after the definition of $\mathfrak{M}$, \eqref{eq:LVD} implies that the operator $\mathfrak{M}$ is bounded on $L^2$, so using H\"older, (3) from the Calder\'on-Zygmund decomposition, and  (\ref{lem:max}) we get that

\begin{eqnarray*}
\left\|  \chi_{3B}\sum_{i\in I} \int_{M}\frac{\mathrm{e}^{-\frac{\varrho^2(\cdot,y)}{Ct_i}}}{V(y,\sqrt{t_i})}\chi_{B_i}(y)\, \Id\mu(y) \right\|_2^2 &\leq & C ||\mathfrak{M}v||_2^2 \sum_{i\in I}\mu(B_i)\\
&\leq & \frac{C}{\lambda^2}||v||_2^2||u||_1.
\end{eqnarray*}
Dividing by $||v||_2^2$ and taking the sup over all non-zero $v$, we obtain

$$\left\|  \chi_{3B}\sum_{i\in I} \int_{M}\frac{\mathrm{e}^{-\frac{\varrho^2(\cdot,y)}{Ct_i}}}{V(y,\sqrt{t_i})}\chi_{B_i}(y)\, \Id\mu(y) \right\|_2^2 \leq \frac{C}{\lambda}||u||_1,$$
which proves \eqref{eq:Riesz22}, hence \eqref{eq:Riesz2}. 

It thus remains to prove \eqref{eq:Riesz3} and \eqref{eq:Riesz4}. It relies on the following lemma:

\begin{Lemma}\label{lem:grad-Riesz} Assume $\left\|\mathrm{Riem}\right\|_{\infty}\leq A<\infty$. Then there is a constant $C=C(A,m)>0$, such that for every $t>0$, $s>0$ and $y\in M$,
$$
\int_{\{\varrho(\cdot,y)\geq \sqrt{t}\}} |(\Id_j+\Id^\dagger_{j-1}) \mathrm{e}^{-s\vec{\Delta}_j}(x,y)|\, \Id\mu(x)\leq Cs^{-1/2}\mathrm{e}^{-\frac{t}{Cs}}\mathrm{e}^{Cs}.
$$

Assume 
$$
\max(\left\|\mathrm{Riem}\right\|_{\infty},\left\|\nabla\mathrm{Riem}\right\|_{\infty})\leq A<\infty.
$$
 Then there is a constant $C=C(A,m)>0$, such that for every $t>0$, $s>0$ and $y\in M$,
$$\int_{\{\varrho(\cdot,y)\geq \sqrt{t}\}} |\nabla \mathrm{e}^{-s\vec{\Delta}_j}(x,y)|\,\Id\mu(x)\leq Cs^{-1/2}\mathrm{e}^{-\frac{t}{Cs}}\mathrm{e}^{Cs}.
$$

\end{Lemma}

\begin{proof}
For the integral involving $\nabla \mathrm{e}^{-s\vec{\Delta}_j}(x,y)$, it is an immediate consequence of Corollary \ref{cor:gradient}, part II with the choice $p=1$. The proof for the second integral follows along the same lines, using \eqref{eq:dj} and \eqref{eq:d*j} instead of \eqref{eq:CUE} for the proof of the weighted estimate analogous to  Corollary \ref{cor:gradient}, part II.
\end{proof}
The estimates \eqref{eq:Riesz3} and \eqref{eq:Riesz4} follow from Lemma \ref{lem:grad-Riesz}, in a fashion that is identical to the proof of \cite[Theorem 1.2]{cd}, and thus whose details will be omitted. Finally, all four estimates \eqref{eq:Riesz1}, \eqref{eq:Riesz2}, \eqref{eq:Riesz3} and \eqref{eq:Riesz4} are proved. This concludes the proof of Theorem \ref{thm:riesz-bakry} and Corollary \ref{cor:Riesz}.


\section{Proof of Theorem \ref{main2}}\label{beweis3}

We prepare the proof with the following estimate from complex analysis can be found in \cite{CS}:

\begin{Lemma}[Phragmen-Lindelöf's inequality] Let 
$$
f:\{\Re >0\}\longrightarrow \IC
$$
be holomorphic, and assume that there are constants $A,B,\gamma>0$, $b\geq 0$, such that
\begin{align}\label{p1}
&|f(z)|\leq B\>\text{ for all $z\in \{\Re >0\}$},\\
\label{p2}
&|f(t)|\leq A\mathrm{e}^{bt-\f{\gamma}{t} }\>\text{ for all $t>0$.}
\end{align}
Then one has 
\begin{align}
|f(z)|\leq B\mathrm{e}^{-\Re \f{\gamma}{z}  }\>\text{ for all $z\in \{\Re >0\}$}.
\end{align}
\end{Lemma}

\begin{proof}[Proof of Theorem \ref{main2}] Step 1: One has 
$$
\left\|1_F \mathrm{e}^{-t \vec{\Delta}_j}\alpha\right\|_2\leq \mathrm{e}^{-\frac{ \varrho(E,F)^2}{4t}} \left\|1_E\alpha\right\|_{2}.
$$
Proof of step 1: The inequality 
\begin{align*}
\left| \left\langle \mathrm{e}^{-t \vec{\Delta}_j}\alpha_1,\alpha_2\right\rangle\right|\leq \mathrm{e}^{C(A)t} \mathrm{e}^{-\frac{ \varrho(E,F)^2}{4t}} \left\|\alpha_1\right\|_{2}\left\|\alpha_2\right\|_{2}
\end{align*}
valid for all $\alpha_1$ with support in $E$ and $\alpha_2$ with support in $F$ has been proved in \cite{batu2}. If we apply Phragmen-Lindelöf's inequality with
\begin{align*}
f(z)=  \left\langle \mathrm{e}^{-z \vec{\Delta}_j}\alpha_1,\>\alpha_2\right\rangle ,\> b=|a|,\>A=B=\left\|\alpha_1\right\|_2\left\|\alpha_2\right\|_2,\>\gamma=\varrho(E,F)^2/4,
\end{align*}
noting that one may pick $A=\left\|\alpha_1\right\|_2\left\|\alpha_2\right\|_2$ because of $\vec{\Delta}_j\geq 0$ so that $\mathrm{e}^{-z\vec{\Delta}_j}$ is a contraction, we get the bound
\begin{align} \label{gds4}
\left|\left\langle \mathrm{e}^{-t \vec{\Delta}_j}\alpha_1,\alpha_2 \right\rangle\right|\leq \mathrm{e}^{\f{-\varrho(E,F)^2}{4t}}\left\|\alpha_1\right\|_2\left\|\alpha_2\right\|_2.
\end{align}
The latter inequality is equivalent to the statement of step 1.\vspace{1mm}

Step 2: One has 
$$
\left\|1_Ft\vec{\Delta}_j  \mathrm{e}^{-t \vec{\Delta}_j}\alpha\right\|_2\leq C\mathrm{e}^{-\frac{ \varrho(E,F)^2}{6t}} \left\|1_E\alpha\right\|_{2},
$$
where $C<\infty$ is a universal constant.\\
Proof of step 2: The asserted estimate is equivalent to 
\begin{align} \label{gds5}
\left|\left\langle t\vec{\Delta}_j\mathrm{e}^{-t \vec{\Delta}_j}\alpha_1,\alpha_2 \right\rangle\right|\leq C\mathrm{e}^{ -\f{\varrho(U_1,U_2)^2}{6t}}\left\|\alpha_1\right\|_2\left\|\alpha_2\right\|_2,
\end{align}
where $\alpha_1\in\Gamma_{L^2}(M,\Lambda^jT^*M)$ is supported in $E$ and $\alpha_2\in\Gamma_{L^2}(M,\Lambda^jT^*M)$ is supported in $F$. To see (\ref{gds5}), we first note that by applying the Phragmen-Lindelöf estimate to the estimate (\ref{gds4}) we get the bound
\begin{align} \label{gds6}
\left|\left\langle \mathrm{e}^{-z \vec{\Delta}_j}\alpha_1,\alpha_2 \right\rangle\right|\leq \mathrm{e}^{-\varrho(U_1,U_2)^2\Re\f{1}{4z}}\left\|\alpha_1\right\|_2\left\|\alpha_2\right\|_2,
\end{align}
valid for all $z$ with $\Re z>0$. By Cauchy's integral formula we have
\begin{align}
\left\langle \vec{\Delta}_j\mathrm{e}^{-t \vec{\Delta}_j}\alpha_1,\alpha_2 \right\rangle= -\frac{\Id}{\Id t} \left\langle \mathrm{e}^{-t \vec{\Delta}_j}\alpha_1,\alpha_2 \right\rangle = -\frac{1}{2\pi i} \int_{z: |z-t|= t/2} \frac{ \left\langle \mathrm{e}^{-z \vec{\Delta}_j}\alpha_1,\alpha_2 \right\rangle \Id z}{(z-t)^2}.
\end{align}
Since by (\ref{gds6}) we have
\begin{align*}
\left|\int_{z: |z-t|= t/2} \frac{ \left\langle \mathrm{e}^{-z \vec{\Delta}_j}\alpha_1,\alpha_2 \right\rangle \Id z}{(z-t)^2}\right|&\leq (2\pi)^{-1}  \pi t \sup_{z: |z-t|=t/2} \left|\frac{ \left\langle \mathrm{e}^{-z \vec{\Delta}_j}\alpha_1,\alpha_2 \right\rangle }{(z-t)^2}\right|\\
&\leq (1/2) t \left\|\alpha_1\right\|_2\left\|\alpha_2\right\|_2 \mathrm{e}^{-\f{\varrho(U_1,U_2)^2}{4(t+t/2)}}(t/2)^{-2},
\end{align*}
this proves step 2. \vspace{1mm}

Step 3: One has 
$$
\left\|1_F\sqrt{t} \nabla \mathrm{e}^{-t \vec{\Delta}_j}\alpha\right\|_2\leq C_1(A)\mathrm{e}^{-\frac{C_2(A) \varrho(E,F)^2}{t}} \left\|1_E\alpha\right\|_{2}.
$$

Proof of step 3: Pick $\phi\in C^{\infty}_c(M)$. Then we have
\begin{align*}
&\left\| \sqrt{t}\phi\nabla \mathrm{e}^{-t \vec{\Delta}_j}\alpha\right\|^2_2 \\
&=\left\langle t\nabla^{\dagger} ( \phi^2\nabla \mathrm{e}^{-t \vec{\Delta}_j}\alpha), \mathrm{e}^{-t \vec{\Delta}_j}\alpha\right\rangle\\
&=2\left\langle t\phi \nabla_{\Id\phi} \mathrm{e}^{-t \vec{\Delta}_j}\alpha, \mathrm{e}^{-t \vec{\Delta}_j}\alpha\right\rangle + \left\langle t\phi^2\nabla^{\dagger} \nabla \mathrm{e}^{-t \vec{\Delta}_j}\alpha, \mathrm{e}^{-t \vec{\Delta}_j}\alpha\right\rangle\\
&=2\left\langle t\phi \nabla \mathrm{e}^{-t \vec{\Delta}_j}\alpha, \Id\phi\otimes \mathrm{e}^{-t \vec{\Delta}_j}\alpha\right\rangle + \left\langle t\phi^2 \vec{\Delta}_j \mathrm{e}^{-t \vec{\Delta}_j}\alpha, \mathrm{e}^{-t \vec{\Delta}_j}\alpha\right\rangle-\left\langle t\phi^2V_j \mathrm{e}^{-t \vec{\Delta}_j}\alpha, \mathrm{e}^{-t \vec{\Delta}_j}\alpha\right\rangle\\
&\leq 2 \left\|\sqrt{t}\phi \nabla \mathrm{e}^{-t \vec{\Delta}_j}\alpha \right\|_2 \sqrt{t}\left\|\Id\phi\otimes \mathrm{e}^{-t \vec{\Delta}_j}\alpha\right\|_2+  \left\|t\phi \vec{\Delta}_j \mathrm{e}^{-t \vec{\Delta}_j}\alpha\right\|_2\left\| \phi \mathrm{e}^{-t \vec{\Delta}_j}\alpha\right\|_2+ A^2t\left\| \phi \mathrm{e}^{-t \vec{\Delta}_j}\alpha \right\|_2^2\\
&\leq  (1/2)\left\|\sqrt{t}\phi \nabla \mathrm{e}^{-t \vec{\Delta}_j}\alpha \right\|_2^2+ 4 t\left\| \Id\phi\otimes \mathrm{e}^{-t \vec{\Delta}_j}\alpha\right\|^2_2+  \left\|\phi t\vec{\Delta}_j \mathrm{e}^{-t \vec{\Delta}_j}\alpha\right\|_2\left\| \phi \mathrm{e}^{-t \vec{\Delta}_j}\alpha\right\|_2+ A^2t\left\| \phi \mathrm{e}^{-t \vec{\Delta}_j}\alpha \right\|_2^2,
\end{align*}
and so 
\begin{align*}
\left\| \sqrt{t}\phi\nabla \mathrm{e}^{-t \vec{\Delta}_j}\alpha\right\|^2_2\leq ct\left\| \Id\phi\otimes \mathrm{e}^{-t \vec{\Delta}_j}\alpha\right\|^2_2+ c \left\|t\phi \vec{\Delta}_j \mathrm{e}^{-t \vec{\Delta}_j}\alpha\right\|_2\left\| \phi \mathrm{e}^{-t \vec{\Delta}_j}\alpha\right\|_2+ cA^2t\left\|t \phi \mathrm{e}^{-t \vec{\Delta}_j}\alpha \right\|_2^2,
\end{align*}
Assume now that 
$$
0\leq \phi\leq 1,\quad \phi|_F=1,\quad \left\|\Id\phi\right\|_{\infty}\leq  1,\quad \mathrm{supp}(\phi)\subset F':=\{x: \varrho(x,F)\leq \varrho(E,F)/3\}.
$$
Then we have
\begin{align*}
&\left\| 1_F\sqrt{t}\nabla \mathrm{e}^{-t \vec{\Delta}_j}\alpha\right\|^2_2 \leq \left\| \sqrt{t}\phi\nabla \mathrm{e}^{-t \vec{\Delta}_j}\alpha\right\|^2_2 \\
 &\leq ct\left\| 1_{F'} \mathrm{e}^{-t \vec{\Delta}_j}\alpha\right\|^2_2+  c\left\|1_{F'} t\vec{\Delta}_j \mathrm{e}^{-t \vec{\Delta}_j}\alpha\right\|_2\left\| 1_{F'} \mathrm{e}^{-t \vec{\Delta}_j}\alpha\right\|_2+ ctA^2\left\|  1_{F'} \mathrm{e}^{-t \vec{\Delta}_j}\alpha \right\|_2^2.
\end{align*}
Using step 1 and step 2 and
$$
\varrho(E,F')\geq \frac{2}{3}\varrho(E,F)
$$
we get
\begin{align*}
&\left\| 1_F\sqrt{t}\nabla \mathrm{e}^{-t \vec{\Delta}_j}\alpha\right\|_2\leq c_1(1+\sqrt{t}A)\mathrm{e}^{-\frac{c_2 \varrho(E,F)^2}{t}} \left\|1_E\alpha\right\|_{2}.
\end{align*}
\end{proof}

\section{Proof of Theorem \ref{main3}}\label{beweis4}

The assumption is equivalent to the following operator norm bound 
\begin{equation}\label{DG}
||1_F\sqrt{t}\nabla \mathrm{e}^{-t\vec{\Delta}_j}1_E||_{2,2}\lesssim \mathrm{e}^{-\frac{\varrho(E,F)^2}{Ct}}
\end{equation}
for disjoint Borel subsets $E, F\subset M $ with compact closure and every $t>0$. Let $A_t=\sqrt{t}\nabla \mathrm{e}^{-t\vec{\Delta}_j}$. Fix an arbitrary $E$ as above, fix $t>0$, and define
$$
F_0:= \left\{x\in M\,:\,\varrho(x,E)^2/t\leq 1\right\}
$$
and for $n\geq1$,
$$
F_n:=\left\{x\in M\,:\,2^{n-1}<\varrho(x,E)^2/t\leq 2^n \right\}.
$$
Clearly,
$$
M=\bigsqcup_{n=0}^\infty F_n,
$$
thus, by \eqref{DG},
$$
||A_t1_E||_{2,2}\leq \sum_{n=0}^\infty ||1_{F_n}A_t1_E||_{2,2}\lesssim \sum_{n=0}^\infty \mathrm{e}^{-c2^{-n}}<+\infty.
$$
So, we have that the operator $A_t1_E$ is bounded in $\Gamma_{L^2}(M,\Lambda^jT^*M)$, uniformly with respect to $t>0$ and the set $E$. Taking an exhaustion of $M$ by compacts $E_n\nearrow M$, we obtain that $A_t$ is bounded in $\Gamma_{L^2}(M,\Lambda^jT^*M)$, uniformly in $t>0$. Let $\alpha\in\Gamma_{L^2}(M,\Lambda^jT^*M)$. We compute, using $\vec{\Delta}_j=\nabla^{\dagger}\nabla + V_j$,
$$\begin{array}{rcl}
0\leq ||\sqrt{t}\nabla \mathrm{e}^{-t\vec{\Delta}_j}\alpha||^2_{2}&=&\left\langle t\nabla^{\dagger}\nabla \mathrm{e}^{-t\vec{\Delta}_j}\alpha,\mathrm{e}^{-t\vec{\Delta}_j}\alpha\right\rangle\\
&=&  \left\langle t\vec{\Delta}_j \mathrm{e}^{-t\vec{\Delta}_j}\alpha,\mathrm{e}^{-t\vec{\Delta}_j}\alpha\right\rangle-t\int_M(V_j\mathrm{e}^{-t\vec{\Delta}_j}\alpha,\mathrm{e}^{-t\vec{\Delta}_j}\alpha)\Id\mu.
\end{array}$$
The first term on the right is bounded, uniformly in $t>0$ by the spectral theorem. Also, the left hand side is non-negative and bounded. Let us take $\alpha\in \mathrm{Ker}_{L^2}(\vec{\Delta}_j)$, then $\mathrm{e}^{-t\vec{\Delta}_j}\alpha=\alpha$, and we get that
$$
t\int_M(V_j\alpha,\alpha)\Id\mu
$$
is non-positive and bounded, which can happen only if
$$
\int_M(V_j\alpha,\alpha)\Id\mu=0.
$$
Since $\vec{\Delta}_j\alpha=0$, we get by the Bochner formula that
$$
0=\left\langle \vec{\Delta}_j\alpha,\alpha\right\rangle=\left\langle \nabla^{\dagger}\nabla\alpha,\alpha\right\rangle+\int_M(V_j\alpha,\alpha)\Id\mu= \left\langle \nabla^{\dagger}\nabla\alpha,\alpha\right\rangle.$$
Thus,
$$ 
\left\langle \nabla^{\dagger}\nabla\alpha,\alpha\right\rangle=0.
$$
Since $V_j$ is assumed to be bounded and since $\alpha\in\dom(\Delta_j)$, it follows that $\alpha\in\dom(\nabla^{\dagger}\nabla)$, and by integration by parts
$$
||\nabla \alpha||_2^2=\left\langle \nabla^{\dagger}\nabla\alpha,\alpha\right\rangle=0.
$$
Thus, $\nabla\alpha=0$, and $\alpha$ is parallel. Since parallel transport with respect to the Levi-Civit\`a connection is an isometry, it follows that
$$
|\alpha(x)|=|\alpha(y)|\quad\text{for all $x,y\in M$},
$$
and since $|\alpha|$ is in $L^2(M)$ and $M$ is non-compact, we conclude that $\alpha\equiv0$.

\appendix

\section{Proof of Proposition \ref{pro:gradUE}}  

Recall that under the assumptions of Proposition \ref{pro:gradUE}, which we assume everywhere in this section, \eqref{eq:vUE} holds. We divide the proof of Proposition \ref{pro:gradUE} into a sequence of lemmas. For $\alpha>0$, $V_{\sqrt{t}}^\alpha$ will denote the operator of multiplication by the function $x\mapsto V(x,\sqrt{t})^\alpha$. Then, we consider two inequalities for $1\leq r\leq s\leq+\infty $; first look at the following $L^r\to L^s$ estimate:

\begin{equation}\label{eq:gradUE_rs}\tag{$\mathrm{G}_{r,s}$}
\sup_{t>0}\mathrm{e}^{-Ct}||V_{\sqrt{t}}^{\frac{1}{r}-\frac{1}{s}}\sqrt{t}(\Id_j+\Id^\dagger_{j-1})\mathrm{e}^{-t\vec{\Delta}_j}||_{r, s}<+\infty
\end{equation}
and next look at its off-diagonal counterpart: for all $x,y\in M,$ and $t>0$,

\begin{equation}\label{eq:gradUE_off}\tag{$\mathrm{G}^{\mathrm{off}}_{r,s}$}
\mathrm{e}^{-Ct}V(x,\sqrt{t})^{\frac{1}{r}-\frac{1}{s}}||\chi_{B(x,\sqrt{t})}\sqrt{t}(\Id_j+\Id^\dagger_{j-1})\mathrm{e}^{-t\vec{\Delta}_j}\chi_{B(y,\sqrt{t})}||_{r, s}\leq C\mathrm{e}^{-\frac{\varrho^2(x,y)}{Ct}},
\end{equation}
where, for a Borel set $A$, $\chi_A$ denotes the indicator function of $A$. Note that if $w\in B(x,\sqrt{t})$, then according to the volume comparison \eqref{eq:VC}, one has

$$\frac{V(x,\sqrt{t})}{V(w,\sqrt{t})}\leq C\mathrm{e}^{Ct},$$
therefore the estimate in \eqref{eq:gradUE_off} follows from the one in \eqref{eq:gradUE_rs} (with a different constant $C$) provided $t\geq C \varrho^2(x,y)$. For small times, one has the following lemma:

\begin{Lemma}\label{lem:off-diag}

Assume that for some $s\geq 2$ (resp., $r\leq 2$), {\em ($\mathrm{G}_{2,s}$)} (resp. {\em ($\mathrm{G}_{r,2}$)}) holds. Then {\em ($\mathrm{G}^{\mathrm{off}}_{2,s}$)} (resp. {\em ($\mathrm{G}^{\mathrm{off}}_{r,2}$)}) holds.

\end{Lemma}

\begin{proof}

The proof of the case $r\leq 2$ follows by duality from the case $s\geq2$, so we only prove the latter. Let $r=\varrho(x,y)$ and $\alpha=\frac{1}{2}-\frac{1}{s}$. Let $f_1\in \Gamma_{L^2}(M,\Lambda^jT^*M)$ and 
$$
f_2=\omega_2+\eta_2\in \Gamma_{L^2\cap L^s}(M,\Lambda^{j-1}T^*M\oplus \Lambda^{j+1}T^*M)
$$
be smooth, such that the support of $f_1$ (resp. $f_2$) is included in $B(y,\sqrt{t})$ (resp. $B(x,\sqrt{t})$). For $z\in H:=\{\mathfrak{R}>0\}$, set $t=\mathfrak{R}(z)$. Let
$$F(z)=\langle f_2,(\Id_j+\Id^\dagger_{j-1})\mathrm{e}^{-z\vec{\Delta}_j}f_1\rangle.$$
Writing
$$F(z)=\langle \Id^\dagger_{j-1} \eta_2+\Id_j\omega_2,\mathrm{e}^{-z\vec{\Delta}_j}f_1\rangle,$$
the spectral theorem implies that for all $z\in H$,
\begin{equation}\label{eq:unif}
|F(z)|\leq ||f_1||_2||\Id^\dagger_{j-1} \eta_2+\Id_j\omega_2||_2<+\infty,
\end{equation}
hence $F$ is uniformly bounded on $H$. According to the Davies-Gaffney estimate for $\sqrt{t}(\Id_j+\Id^\dagger_{j-1})\mathrm{e}^{-z\vec{\Delta}_j}$ (\cite[Lemma 3.8]{AMR}), for all $t>0$, it holds:
$$|F(t)|\leq \frac{C}{\sqrt{t}}\mathrm{e}^{-\frac{r^2}{Ct}}||f_1||_2||f_1||_2,$$
hence, if one lets $\gamma=\frac{r^2}{2C}$, one has for all $0\leq t\leq \gamma$,
\begin{equation}\label{eq:estDG}
|F(t)|\leq \frac{C}{r}\mathrm{e}^{-\frac{\gamma}{t}}||f_1||_2||f_2||_2.
\end{equation}
As explained before the statement of Lemma \ref{lem:off-diag}, one can limit ourselves to prove ($\mathrm{G}^{\mathrm{off}}_{2,s}$) for all $t\leq \gamma$, hence in the rest of the proof we will assume that $t\leq \gamma$. Next, write
$$
F(z)=V(x,r)^{-\alpha}\left\langle \left(\frac{V(x,r)}{V(\cdot,\sqrt{t})}\right)^\alpha f_2,V_{\sqrt{t}}^\alpha(\Id_j+\Id^\dagger_{j-1})\mathrm{e}^{-z\vec{\Delta}_j}f_1\right\rangle.
$$
Writing $z=t+iu$, $\mathrm{e}^{-z\vec{\Delta}_j}=\mathrm{e}^{-t\vec{\Delta}_j}\mathrm{e}^{-iu\vec{\Delta}_j}$, and using that $||\mathrm{e}^{-iu\vec{\Delta}_j}||_{2,2}\leq 1$ (by self-adjointness and the spectral theorem) and ($\mathrm{G}_{2,s}$), one has
$$||V_{\sqrt{t}}^\alpha(\Id_j+\Id^\dagger_{j-1})\mathrm{e}^{-z\vec{\Delta}_j}f_1||_s\leq \frac{C}{\sqrt{t}}\mathrm{e}^{Ct}||f_1||_2.$$
For any $w$ in the support of $f_2$ (which is included inside the ball $B(x,r)$), one has by volume comparison \eqref{eq:VC},
$$\frac{V(x,r)}{V(w,\sqrt{t})}\leq C\left(\frac{r}{\sqrt{t}}\right)^{n}\mathrm{e}^{Cr},$$
therefore we obtain
\begin{eqnarray}\label{eq:est_C}
|F(z)|&\leq & C\frac{V(x,r)^{-\alpha}}{\sqrt{t}}\left(\frac{t}{\gamma}\right)^{-n/2}\mathrm{e}^{Cr}||f_1||_2||f_2||_{s'},\\
&\leq & C\frac{V(x,r)^{-\alpha}}{r}\left(\frac{t}{\gamma}\right)^{-\frac{n+1}{2}}\mathrm{e}^{Cr}||f_1||_2||f_2||_{s'},
\end{eqnarray}
where $s'$ is the conjugate to $s$. According to \cite[Proposition 2.3]{CS}, \eqref{eq:unif}, \eqref{eq:estDG} and \eqref{eq:est_C} imply: for all $t\leq \gamma$, 
\begin{eqnarray*}
|F(t)| &\leq & C\frac{V(x,r)^{-\alpha}}{r}\mathrm{e}^{Cr}\mathrm{e}^{-\frac{r^2}{Ct}}||f_1||_2||f_2||_{s'}\\
&\leq & C\frac{V(x,r)^{-\alpha}}{\sqrt{t}}\mathrm{e}^{Cr}\mathrm{e}^{-\frac{r^2}{Ct}}||f_1||_2||f_2||_{s'}
\end{eqnarray*}
Making use of the elementary inequality
\begin{equation}\label{eq:elem}
\mathrm{e}^{Cr}\leq \mathrm{e}^{\frac{\epsilon r^2}{2t}}\mathrm{e}^{\frac{Ct}{2\epsilon}},
\end{equation}
with a choice of $\epsilon$ small enough, we obtain, for all $t\leq \gamma$,
$$|F(t)|\leq C\frac{V(x,r)^{-\alpha}}{\sqrt{t}}\mathrm{e}^{Ct}\mathrm{e}^{-\frac{r^2}{Ct}}||f_1||_2||f_2||_{s'}.$$
A density argument then yields that the same inequality holds for  
$$
\text{$f_2$ merely in}\quad\Gamma_{L^s}(M,\Lambda^{j-1}T^*M\oplus \Lambda^{j+1}T^*M),
$$
and taking the supremum over $f_1\in \Gamma_{L^2}(M,\Lambda^jT^*M)$ one concludes that ($\mathrm{G}^{\mathrm{off}}_{2,s}$) holds.
\end{proof}

Next, we state an $L^2$-estimate, which actually holds without any further assumptions on the geometry, and which will be useful for the proof of ($\mathrm{G}_{2,s}$):

\begin{Lemma}\label{lem:est_L2}

The following $L^2$-estimate holds:
$$
\sup_{t>0}||\sqrt{t}(\Id_j+\Id^\dagger_{j-1})\mathrm{e}^{-t\vec{\Delta}_j}||_{2,2}\leq A,
$$
where $A<\infty$ is a universal constant.
\end{Lemma}

\begin{proof} 

Let $f\in \Gamma_{L^2}(M,\Lambda T^*M)$. It follows from the fact that $\vec{\Delta}=(\Id+\Id^\dagger)^2$, integration by parts and the spectral theorem that
\begin{eqnarray*}
||\sqrt{t}(\Id+\Id^\dagger)\mathrm{e}^{-t\vec{\Delta}}f||^2_2 &=& \langle t \vec{\Delta}\mathrm{e}^{-t\vec{\Delta}}f,\mathrm{e}^{-t\vec{\Delta}_j}f\rangle \\
& \leq & || t  \vec{\Delta}\mathrm{e}^{-t\vec{\Delta}}f||_2||f||_2\\
&\leq & \left(\sup_{x\geq 0}x\mathrm{e}^{-x}\right)||f||_2^2\\
&\leq & C||f||_2^2,
\end{eqnarray*}
and the result follows.
\end{proof}

We will need the following $L^r\to L^s$ estimates for $\mathrm{e}^{-t\vec{\Delta}_j}$, which extends an analogous result from \cite{BCS} to the locally doubling case:

\begin{Lemma}\label{lem:Lr-Ls}

There is a constant $C$, which only depends on an upper bound of $\left\|\mathrm{Riem}\right\|_{\infty}$, $\left\|\nabla\mathrm{Riem}\right\|_{\infty}$ and on $m$, such that for every $1\leq r\leq s\leq +\infty$, 
\begin{equation}\label{VEV}\tag{$\mathrm{VE}_{r,s}$}
\sup_{t>0}\mathrm{e}^{-Ct}||V_{\sqrt{t}}^{\frac{1}{r}-\frac{1}{s}}\mathrm{e}^{-t\vec{\Delta}_j}||_{r,s}<+\infty.
\end{equation}
\end{Lemma}

\begin{proof}

An interpolation argument (see \cite[Proposition 2.1.5]{BCS}) implies that it is enough to prove ($\mathrm{VE}_{1,1}$), ($\mathrm{VE}_{\infty,\infty}$) and ($\mathrm{VE}_{1,\infty}$). However, the latter is equivalent to the pointwise bound on the kernel $\mathrm{e}^{-t\vec{\Delta}_j}(x,y)$:
$$|\mathrm{e}^{-t\vec{\Delta}_j}(x,y)|\leq \frac{C\mathrm{e}^{Ct}}{V(x,\sqrt{t})},$$
which holds by virtue of \eqref{eq:vUE}. Thus, ($\mathrm{VE}_{1,\infty}$) holds. Concerning ($\mathrm{VE}_{1,1}$) and ($\mathrm{VE}_{\infty,\infty}$), they follow from \eqref{eq:Lpaction} with $p=1$ and $p=\infty$ respectively. The proof is complete.

\end{proof}
The next result concerns the validity of ($\mathrm{G}_{2,s}$) in the context of Proposition \ref{pro:gradUE}:

\begin{Lemma}\label{lem:G2s}

There exists a constant $C$, which only depends on an upper bound of $\left\|\mathrm{Riem}\right\|_{\infty}$, $\left\|\nabla\mathrm{Riem}\right\|_{\infty}$ and on $m$, such that for all $s\in [2,+\infty]$, {\em ($\mathrm{G}_{2,s}$)} holds.

\end{Lemma}

\begin{proof}

We use the commutation rules $\Id_j\mathrm{e}^{-t\vec{\Delta}_j}=\mathrm{e}^{-t\vec{\Delta}_{j+1}}\Id_j$ as well as $\Id^\dagger_{j-1} \mathrm{e}^{-t\vec{\Delta}_j}=\mathrm{e}^{-t\vec{\Delta}_{j-1}}\Id^\dagger_{j-1}$. Write
$$
\sqrt{t}(\Id_j+\Id^\dagger_{j-1})\mathrm{e}^{-t\vec{\Delta}_j}=\mathrm{e}^{-\frac{t}{2}\vec{\Delta}_{j+1}}\left(\sqrt{t}\Id_{j-1}\mathrm{e}^{-\frac{t}{2}\vec{\Delta}_{j-1}}\right)+\mathrm{e}^{-\frac{t}{2}\vec{\Delta}_{j-1}}\left(\sqrt{t}\Id^\dagger_j \mathrm{e}^{-\frac{t}{2}\vec{\Delta}_{j+1}}\right)
$$
Lemma \ref{lem:est_L2} implies that the operators $\sqrt{t}\Id^\dagger_j \mathrm{e}^{-\frac{t}{2}\vec{\Delta}_{j+1}}$ and $\sqrt{t}\Id_{j-1}\mathrm{e}^{-\frac{t}{2}\vec{\Delta}_{j-1}}$ are uniformly bounded on $L^2$. On the other hand, Lemma \ref{lem:Lr-Ls} implies that there is a constant $C$ such that for every $s\in [2,+\infty]$, and for $k\in \{j-1,j+1\}$,
$$
\sup_{t>0}\mathrm{e}^{-Ct}||V_{\sqrt{t}}^{\frac{1}{2}-\frac{1}{s}}\mathrm{e}^{-t\vec{\Delta}_k}||_{2,s}<+\infty.
$$
By composition, we get
$$
\sup_{t>0}\mathrm{e}^{-Ct}||V_{\sqrt{\frac{t}{2}}}^{\frac{1}{2}-\frac{1}{s}}\sqrt{t}(\Id_j+\Id^\dagger_{j-1})\mathrm{e}^{-t\vec{\Delta}_j}||_{2,s}<+\infty.
$$
Local doubling \eqref{eq:LVD} implies that there is a constant $C$ such that for all $s\geq 2$, all $x\in M$ and all $t>0$,
\begin{equation}\label{eq:localvol}
\left(\frac{V(x,\sqrt{t})}{V(x,\sqrt{\frac{t}{2}})}\right)^{\frac{1}{2}-\frac{1}{s}}\leq C\mathrm{e}^{Ct},
\end{equation}
hence one gets ($\mathrm{G}_{2,s}$).

\end{proof}
The last lemma that will be needed for the proof of Proposition \ref{pro:gradUE} is the following composition lemma for Gaussian kernels under local doubling:

\begin{Lemma}\label{lem:compose}

Let $(T_t)_{t>0}$ be a family of operators acting between sections of Riemannian vector bundles over $M$, such that for all $x,y,z\in M$ and $t>0$,
\begin{equation}\label{eq:2-inf}
||\chi_{B(x,\sqrt{t})}T_t\chi_{B(y,\sqrt{t})}||_{2,\infty}\leq \frac{C\mathrm{e}^{Ct}}{V(x,\sqrt{t})^{1/2}}\mathrm{e}^{-\frac{\varrho^2(x,y)}{Ct}},
\end{equation}
and
\begin{equation}\label{1-2}
||\chi_{B(y,\sqrt{t})}T_t\chi_{B(z,\sqrt{t})}||_{1,2}\leq \frac{C\mathrm{e}^{Ct}}{V(z,\sqrt{t})^{1/2}}\mathrm{e}^{-\frac{\varrho^2(y,z)}{Ct}}.
\end{equation}
Then there is a constant $C>0$, which only depends on the local doubling constant, such that for all $x,y\in M$ and $t>0$ such that $\varrho(x,z)\geq \sqrt{t}$,
\begin{equation}\label{1-inf}
||\chi_{B(x,\sqrt{t})}(T_{t}\circ T_t)\chi_{B(z,\sqrt{t})}||_{1,\infty}\leq \frac{C\mathrm{e}^{Ct}}{V(z,\sqrt{t})^{1/2}V(x,\sqrt{t})^{1/2}}\mathrm{e}^{-\frac{\varrho^2(x,z)}{Ct}}.
\end{equation}

\end{Lemma}
The proof of Lemma \ref{lem:compose} is a bit long, so we postpone it to Appendix B. Finally, one can conclude the proof of Proposition \ref{pro:gradUE}. According to Lemma \ref{lem:G2s}, ($\mathrm{G}_{2,\infty}$), hence (Lemma \ref{lem:off-diag}) ($\mathrm{G}^{\mathrm{off}}_{2,\infty}$), holds. By duality, one gets
\begin{equation}\label{eq:dual}
\mathrm{e}^{-Ct}V(y,\sqrt{t})^{1/2}||\chi_{B(x,\sqrt{t})}\sqrt{t}(\Id_j+\Id^\dagger_{j-1})\mathrm{e}^{-t\vec{\Delta}_j}\chi_{B(y,\sqrt{t})}||_{1,2}\leq C\mathrm{e}^{-\frac{\varrho^2(x,y)}{Ct}},
\end{equation}
as well. By the composition lemma (Lemma \ref{lem:compose}) applied to $T_t=\sqrt{t} (\Id_j+\Id^\dagger_{j-1})\mathrm{e}^{-t\vec{\Delta}_j}$ and local volume doubling \eqref{eq:LVD}, one gets for all $x,y\in M$ and $t>0$ such that $\varrho^2(x,y)\geq \frac{t}{2}$,
$$
||\chi_{B(x,\sqrt{t})}\sqrt{t}(\Id_j+\Id^\dagger_{j-1})\mathrm{e}^{-t\vec{\Delta}_j}\chi_{B(y,\sqrt{t})}||_{1,\infty}\leq \frac{C}{V(x,\sqrt{t})^{1/2}V(y,\sqrt{t})^{1/2}}\mathrm{e}^{-\frac{\varrho^2(x,y)}{Ct}}.
$$
This estimate also holds true for $\varrho^2(x,y)\leq \frac{t}{2}$ (the exponential term becoming a constant), as a consequence of the composition of the estimate ($\mathrm{G}_{2,\infty}$) and its dual:
$$
\sup_{t>0} \mathrm{e}^{-Ct}V(y,\sqrt{t})^{1/2}||\sqrt{t}(\Id_j+\Id^\dagger_{j-1})\mathrm{e}^{-t\vec{\Delta}_j}\chi_{B(y,\sqrt{t})}||_{1,2}<+\infty.
$$
Hence the kernel of $\chi_{B(x,\sqrt{t})}\sqrt{t}(\Id_j+\Id^\dagger_{j-1})\mathrm{e}^{-t\vec{\Delta}_j}\chi_{B(y,\sqrt{t})}$ is pointwise bounded by the right-hand side of the above equation. Evaluating the kernel at $(x,y)$, we obtain:
\begin{equation}\label{eq:L1Linf}
||\sqrt{t}(\Id_j+\Id^\dagger_{j-1})\mathrm{e}^{-t\vec{\Delta}_j}(x,y)||\leq \frac{C}{V(x,\sqrt{t})^{1/2}V(y,\sqrt{t})^{1/2}}\mathrm{e}^{-\frac{\varrho^2(x,y)}{Ct}}.
\end{equation}
Let $r=\varrho(x,y)$, then by local doubling \eqref{eq:LVD} and \eqref{eq:elem}, one has for every $\epsilon>0$,
\begin{eqnarray*}
\frac{V(x,\sqrt{t})}{V(y,\sqrt{t})} &\leq & \frac{V(y,r+\sqrt{t})}{V(y,\sqrt{t})}\\
&\leq & \left(1+\frac{r}{\sqrt{t}}\right)^m\mathrm{e}^{C(r+\sqrt{t})}\\ 
&\leq & \left(1+\frac{r}{\sqrt{t}}\right)^m \mathrm{e}^{C\epsilon^{-1}t}\mathrm{e}^{\epsilon\frac{r^2}{t}}\\
&\leq & C \mathrm{e}^{C\epsilon^{-1}t}\mathrm{e}^{\epsilon\frac{r^2}{t}}.
\end{eqnarray*}
Taking $\epsilon$ small enough and plugging this inequality into \eqref{eq:L1Linf}, we arrive to

$$||\sqrt{t}(\Id_j+\Id^\dagger_{j-1})\mathrm{e}^{-t\vec{\Delta}_j}(x,y)||\leq \frac{C}{V(x,\sqrt{t})}\mathrm{e}^{-\frac{\varrho^2(x,y)}{Ct}},$$
This implies \eqref{eq:dj} and \eqref{eq:d*j}, and it concludes the proof of Proposition \ref{pro:gradUE}.

\section{Proof of Lemma \ref{lem:compose}}

In all the proofs, we write $r:=\varrho(x,z)$. Let $t>0$, and let $(y_i)_{i\in\mathbb{N}}$ be a maximal $\sqrt{t}$-separated set of $M$: this means that for every $i\neq j$, $\varrho(y_i,y_j)\geq \sqrt{t}$, and the set of points $(y_i)_{i\in\mathbb{N}}$ is maximal for this property. Clearly, the balls $\frac{1}{2}B_i:=\frac{1}{2}B(y_i,\sqrt{t})$ are then pairwise disjoint. Moreover, the balls $B_i$, $i\in\mathbb{N}$ cover $M$: indeed, if there is some $y\in M\setminus \cup_{i}B_i$, then the set of points $E:=\{y\}\cup (y_i)_{i\in\mathbb{N}}$ has the property that for every $y,w\in E$ distinct, $\varrho(y,w)\geq \sqrt{t}$, which contradicts the maximality of the set $(y_i)_{i\in\mathbb{N}}$. To sum up, we thus have a covering of $M$ by balls $(B_i)_{i\in \mathbb{N}}$, such that the balls $\frac{1}{2}B_i$ are pairwise disjoints. 

Next, we claim that this covering is {\em locally finite}, and that actually there is a constant $C>0$ such that, for all $i\in \mathbb{N}$,
$$\mathrm{Card}\{j\in\mathbb{N}\,:\,B_i\cap B_j\neq \emptyset\}\leq C\mathrm{e}^{Ct}.$$
Indeed, this follows from the following estimate, with the choice $\alpha=2$:

\begin{Lemma}\label{lem:card}

There exists a constant $C>0$, which only depends on the local doubling constant, such that for all $x\in M$ and $\alpha\geq \frac{1}{2}$,
$$\mathrm{Card}\{i\in\mathbb{N}\,:\,\varrho(x,y_i)\leq \alpha\sqrt{t}\}\leq C\alpha^m \mathrm{e}^{C\alpha\sqrt{t}}\leq C'\alpha^m \mathrm{e}^{C'\alpha }.$$

\end{Lemma}

\begin{proof}

Consider the balls $\frac{1}{2}B_i$, for $i\in A:=\{j\in\mathbb{N}\,:\,\varrho(x,y_j)\leq \alpha\sqrt{t}$. These balls are disjoint, and included in $B(x,(\alpha+\frac{1}{2})\sqrt{t})\subset B(x,2\alpha\sqrt{t})$. Hence,
$$\sum_{i\in A}\mu(\frac{1}{2}B_i)\leq V(x,2\alpha\sqrt{t}).$$
However, \eqref{eq:LVD}, \eqref{eq:VC2} and the fact that $\alpha\geq \frac{1}{2}$ imply that
\begin{eqnarray*}
\frac{V(x,2\alpha\sqrt{t})}{V(y_i,\frac{\sqrt{t}}{2})} &=& \frac{V(x,2\alpha\sqrt{t})}{V(y_i,2\alpha\sqrt{t})}\frac{V(y_i,2\alpha\sqrt{t})}{V(y_i,\frac{\sqrt{t}}{2})}\\
&\leq & C\left(\frac{\varrho(x,y_i)}{\sqrt{t}}+1\right)^m \mathrm{e}^{C\sqrt{t}} \mathrm{e}^{C\varrho(x,y_i)}\\
&\leq & C\alpha^m \mathrm{e}^{C\sqrt{t}} \mathrm{e}^{C\alpha\sqrt{t}}\\
&\leq & C\alpha^m \mathrm{e}^{C\alpha\sqrt{t}}.
\end{eqnarray*}

\end{proof}
Let $(\varphi_i)_{i\in\mathbb{N}}$ be a partition of unity associated with the covering $(B_i)_{i\in \mathbb{N}}$. Write $\varphi_i=\psi_i^2$ for $\psi\geq 0$. Denote by $M_{\psi_i}$ the operator of multiplication with $\psi_i$; given that $0\leq \psi_i\leq 1$, one has $||M_{\psi_i}||_{2,2}\leq ||\psi_i||_\infty\leq 1$. Then, using $\mathbf{1}=\sum_{i\in\mathbb{N}}\psi^2_i$, the assumption that $\psi_i$ has support in $B_i$, and the off-diagonal estimates for $T_t$, one has
\begin{eqnarray*}
||\chi_{B(x,\sqrt{t})}(T_t\circ T_t)\chi_{B(z,\sqrt{t})}||_{1,\infty} &\leq & \sum_{i\in\mathbb{N}} ||\chi_{B(x,\sqrt{t})}T_t \psi_i||_{2,\infty} ||\psi_i T_t\chi_{B(z,\sqrt{t})}||_{1,2}\\
&\leq & \sum_{i\in\mathbb{N}} ||\chi_{B(x,\sqrt{t})}T_t \chi_{B_i}||_{2,\infty} ||M_{\psi_i}||_{2,2}^2 ||\chi_{B_i} T_t\chi_{B(z,\sqrt{t})}||_{1,2}\\
&\leq & \sum_{i\in\mathbb{N}} ||\chi_{B(x,\sqrt{t})}T_t \chi_{B_i}||_{2,\infty}  ||\chi_{B_i} T_t\chi_{B(z,\sqrt{t})}||_{1,2}\\
&\leq & \frac{C\mathrm{e}^{Ct}}{V(x,\sqrt{t})^{1/2}V(z,\sqrt{t})^{1/2}}\sum_{i\in\mathbb{N}} \exp\left(-\left(\frac{\varrho^2(x,y_i)}{Ct}+\frac{\varrho^2(z,y_i)}{Ct}\right)\right).
\end{eqnarray*}
Hence, the proof of Lemma \ref{lem:compose} follows from the following claim: there exists $C>0$ such that for all $x,z\in M$ and $t>0$ such that $\varrho(x,z)\geq 1$,

\begin{equation}\label{eq:sum}
\sum_{i\in\mathbb{N}} \exp\left(-\left(\frac{\varrho^2(x,y_i)}{Ct}+\frac{\varrho^2(z,y_i)}{Ct}\right)\right)\leq C\mathrm{e}^{Ct}\mathrm{e}^{-\frac{\varrho^2(x,z)}{Ct}}.
\end{equation}
In order to prove \eqref{eq:sum}, we split $\sum_{i\in \mathbb{N}}$ into $\sum_{k=-\infty}^{+\infty}\sum_{i\in A_k}$, where $A_k$ is defined as the set of $i\in\mathbb{N}$ such that 
$$2^{k-1}\leq \frac{\varrho(x,y_i)}{\varrho(x,z)}\leq 2^k.$$
We first bound the sum $\sum_{k=-\infty}^{-1}$: if $k\leq -1$ and $i\in A_k$, then
$$\varrho(x,y_i)\leq \frac{1}{2}\varrho(x,z).$$
Therefore, 
$$\varrho(z,y_i)\geq \frac{1}{2}\varrho(x,z),$$
and consequently
$$\exp\left(-\left(\frac{\varrho^2(x,y_i)}{Ct}+\frac{\varrho^2(z,y_i)}{Ct}\right)\right)\leq \mathrm{e}^{-\frac{\varrho^2(x,z)}{Ct}}.$$
It follows that the sum $\sum_{k=-\infty}^{-1}$ is bounded from above by
$$\mathrm{Card}\{i\,:\,\varrho(x,y_i)\leq \frac{1}{2}\varrho(x,z)\}\mathrm{e}^{-\frac{\varrho^2(x,z)}{Ct}}.$$
Lemma \ref{lem:card} with the choice $\alpha=\frac{1}{2}\frac{\varrho(x,z)}{\sqrt{t}}\geq \frac{1}{2}$ yields
\begin{eqnarray*}
\mathrm{Card}\{i\,:\,\varrho(x,y_i)\leq \frac{1}{2}\varrho(x,z)\}&\leq & C \left(\frac{\varrho(x,z)}{\sqrt{t}}\right)^m \mathrm{e}^{C\varrho(x,z)}\\
&\leq & C\mathrm{e}^{C\epsilon^{-1}t}\mathrm{e}^{\epsilon \frac{\varrho^2(x,z)}{t}},
\end{eqnarray*}
for all $\epsilon>0$, where in the last line we have used \eqref{eq:elem}. Hence, taking $\epsilon$ small enough, one concludes that the sum $\sum_{k=-\infty}^{-1}$ is bounded from above by
$$C\mathrm{e}^{Ct}\mathrm{e}^{-\frac{\varrho^2(x,z)}{Ct}}.$$
We now deal with the sum $\sum_{k=0}^\infty$. For every $k\geq 0$ and $i\in A_k$, one has by definition of $A_k$ that
$$2^{k-1}\leq \frac{\varrho(x,y_i)}{\varrho(x,z)},$$
which implies
$$\exp\left(-\left(\frac{\varrho^2(x,y_i)}{Ct}+\frac{\varrho^2(z,y_i)}{Ct}\right)\right)\leq \mathrm{e}^{-2^k\frac{\varrho^2(x,z)}{2Ct}}.$$
So, 
$$\sum_{i\in A_k}\exp\left(-\left(\frac{\varrho^2(x,y_i)}{Ct}+\frac{\varrho^2(z,y_i)}{Ct}\right)\right)\leq \mathrm{Card}(A_k) \mathrm{e}^{-2^k\frac{\varrho^2(x,z)}{2Ct}}.$$
By definition of $A_k$,
$$A_k\subset \{i\in \mathbb{N}\,:\,\varrho(x,y_i)\leq 2^k\varrho(x,z)\},$$
so that, using Lemma \ref{lem:card} with $\alpha=2^k\frac{\varrho(x,z)}{\sqrt{t}}$, we get
\begin{eqnarray*}
\mathrm{Card}(A_k)&\leq & C\left(\frac{2^k\varrho(x,z)}{\sqrt{t}}\right)^m\mathrm{e}^{2^k\varrho(x,z)}\\
&\leq & C\mathrm{e}^{C\epsilon^{-1}t}\mathrm{e}^{-\epsilon\frac{2^k\varrho^2(x,z)}{t}},
\end{eqnarray*}
for any $\epsilon>0$, where in the last line we have used \eqref{eq:elem}. Taking $\epsilon$ small enough, we arrive to
$$\sum_{i\in A_k}\exp\left(-\left(\frac{\varrho^2(x,y_i)}{Ct}+\frac{\varrho^2(z,y_i)}{Ct}\right)\right)\leq C\mathrm{e}^{Ct}\mathrm{e}^{-\frac{2^k\varrho^2(x,z)}{t}},$$
hence
$$\sum_{k=0}^\infty \sum_{i\in A_k}\exp\left(-\left(\frac{\varrho^2(x,y_i)}{Ct}+\frac{\varrho^2(z,y_i)}{Ct}\right)\right)\leq C\mathrm{e}^{Ct}\sum_{k=0}^\infty \mathrm{e}^{-\frac{2^k\varrho^2(x,z)}{t}}.$$
Next, we use the inequality: for all $s>0$,
\begin{equation}\label{eq:exp_sum}
\sum_{k=0}^\infty \mathrm{e}^{-s2^k}\leq C\frac{\mathrm{e}^{-s}}{1-\mathrm{e}^{-cs}}.
\end{equation}
Indeed, writing
$$\sum_{k=0}^\infty \mathrm{e}^s\mathrm{e}^{-s2^k}=\sum_{k=0}^\infty \mathrm{e}^{-s(2^k -1)},$$
and using that $2^k-1\geq Ck$ for $k\geq 0$, we get
\begin{eqnarray*}
\sum_{k=0}^\infty \mathrm{e}^s\mathrm{e}^{-s2^k} &\leq & \sum_{k=0}^\infty \mathrm{e}^{-Csk}\\
&\leq & \sum_{k=0}^\infty (\mathrm{e}^{-Cs})^k\\
&\leq & \frac{1}{1-\mathrm{e}^{-Cs}},
\end{eqnarray*}
which proves \eqref{eq:exp_sum}. Taking $s=\frac{\varrho^2(x,z)}{t}\geq 1$ (by hypothesis) in \eqref{eq:exp_sum}, we obtain
\begin{eqnarray*}
\sum_{k=0}^\infty \sum_{i\in A_k}\exp\left(-\left(\frac{\varrho^2(x,y_i)}{Ct}+\frac{\varrho^2(z,y_i)}{Ct}\right)\right) &\leq & \frac{C\mathrm{e}^{Ct}}{1-\mathrm{e}^{-C}}\mathrm{e}^{-\frac{\varrho^2(x,z)}{Ct}}\\
&\leq & C\mathrm{e}^{Ct}\mathrm{e}^{-\frac{\varrho^2(x,z)}{Ct}}.
\end{eqnarray*}
This concludes the proof of \eqref{eq:sum}, and the proof of Lemma \ref{lem:compose}.

\section{The localized Calder\'on-Zygmund decomposition}

In this appendix we prove Lemma \ref{lem:CZ}, adapting the classical proof (such as written for instance in \cite{St}) by using some ideas from \cite[Appendix B]{DR}. Let $\lambda>0$ be fixed, such that $\lambda>\frac{C}{\mu(B)}\int_B|u|$, where the value of the constant $C>0$ will be chosen in a moment. Let
$$
\Omega:=\left\{x\in M:\ {\mathfrak M}  u(x)>\lambda\right\},
$$
where we recall that the ``local'' maximal function $\mathfrak{M}$ is defined by

$$\mathfrak{M}u(x)=\sup_{\tilde{B}\ni x\,:\,r(\tilde{B})\leq 8}\frac{1}{\mu(\tilde{B})}\int_{\tilde{B}}|u|.$$
It is easily seen that $\Omega$ is an open subset of $M$ and we set $F:=M\setminus \Omega$. We first claim that $\Omega \subset 2B$. Indeed, if $x\notin 2B$, and $\tilde{B}$ is a ball with radius $r(\tilde{B})\leq 8$ containing $x$ and intersecting the support of $u$ (hence intersecting $B$), then $B\subset 3\tilde{B}$, hence $V(B)\leq cV(\tilde{B})$ by local doubling \eqref{eq:LVD}. Consequently,

$$\frac{1}{V(\tilde{B})}\int_{\tilde{B}}| u|\leq c \frac{1}{V(B)}\int_{B}| u| \le cC^{-1}\lambda,$$
hence, if $C\geq c$, one obtains

$$\frac{1}{V(\tilde{B})}\int_{\tilde{B}}| u|\leq \lambda.$$
Taking the supremum over all balls $\tilde{B}$ with radius $r(\tilde{B})\leq 8$ containing $x$, one gets

$$\mathfrak{M} u(x)\le\lambda,$$
and consequently $x\notin \Omega$. Therefore, we have proved that $\Omega \subset 2B$.

For all $x\in \Omega$, let $r_x:=\frac{1}{10}\varrho(x,M\setminus \Omega)$ and $B_x:=B(x,r_x)$, so that $B_x\subset \Omega$, and $\Omega=\bigcup_{x\in \Omega}B_x$. Since the radii of the balls $B_x$ are uniformly bounded, there exists a denumerable collection of points $(x_i)_{i\ge 1}\in \Omega$ such that the balls $B_{x_i}$ are pairwise disjoint and $\Omega= \bigcup_{i\ge 1} 5B_{x_i}$. For all $i$, write $s_i:=5r_{x_i}\leq 1$ and let $B_i=B(x_i,s_i)$. Notice that $B_i\subset 2B$ for all $i$. Furthermore, the balls $\frac{1}{5}B_i$ being disjoint together with local doubling entail that the covering by balls $B_i$ has the finite intersection property (property (2) of the Calder\'on-Zygmund decomposition). And by construction also, $3B_i\cap F\neq \emptyset$ for every $i$. Let $(\chi_i)_{i\geq 1}$ be a partition of unity of $\Omega$, subordinated to the covering $(B_i)_{i\ge 1}$. Then, define

$$b_i=u\chi_i,$$
so that $b_i$ has support in $B_i$. We also let

$$g=\chi_F u=u-\sum_{i\ge 1}b_i$$
(the above sum in fact contains at every point only a finite number of terms, thanks to the finite intersection property of the covering). The Lebesgue differentiation theorem implies that $|g|\leq \lambda$ a.e. on $F$, proving point (4) of the Calder\'on-Zygmund decomposition. Next, since $r(3B_i)\leq 3\leq 8$, the fact that $3B_i\cap F\neq \emptyset$ implies that 

$$\frac{1}{\mu(3B_i)}\int_{3B_i}|u|\leq \lambda.$$
Local doubling then implies that

$$\frac{1}{\mu(B_i)}\int_{B_i}|u|\leq C\lambda,$$
proving (5) of the Calder\'on-Zygmund decomposition. Finally, using the finite intersection property of the covering (point (2) of the Calder\'on-Zygmund decomposition), one has

\begin{eqnarray*}
\sum_{i\in I}\mu(B_i) & \leq & C\mu(\cup_{i\in I}B_i)\\
&\leq & C\mu(\Omega)\\
&\leq & \frac{C}{\lambda}||u||_1,
\end{eqnarray*}
where in the last line we have used the fact that $\mathfrak{M}$ is weak type $(1,1)$. This proves point (3) of the Calder\'on-Zygmund decomposition, and this concludes the proof of Lemma \ref{lem:CZ}.

\end{document}